\documentclass[preprint]{elsarticle} 

\makeatletter
	\def\ps@pprintTitle{%
 	\let\@oddhead\@empty
	\let\@evenhead\@empty
	\def\@oddfoot{\centerline{\thepage}}%
	\let\@evenfoot\@oddfoot}
\makeatother

\usepackage{etoolbox}
\patchcmd{\MaketitleBox}{\footnotesize\itshape\elsaddress\par\vskip36pt}{\footnotesize\itshape\elsaddress\par\parbox[b][36pt]{\linewidth}{\vfill\hfill\textnormal{\today}\hfill\null\vfill}}{}{}%
\patchcmd{\pprintMaketitle}{\footnotesize\itshape\elsaddress\par\vskip36pt}{\footnotesize\itshape\elsaddress\par\parbox[b][36pt]{\linewidth}{\vfill\hfill\textnormal{\today}\hfill\null\vfill}}{}{}%

\usepackage{tikz}
\usepackage{nomencl, multicol}
\usepackage{longtable}
\usepackage{etoolbox}
\usepackage{tikz-imagelabels}
\usepackage{pgfplots} 
\usepackage{pgfgantt}
\usepackage{pdflscape}
\pgfplotsset{compat=newest} 
\pgfplotsset{plot coordinates/math parser=false} 
\newlength\fwidth
\newlength\fheight

\usepackage[
            colorlinks=true,%
            breaklinks=true,%
            linkcolor=blue,
            urlcolor=blue,%
            citecolor=blue,%
            pdftitle={Robust Design Optimization with Limited Data for a Char Combustion Process}, 
            pdfkeywords={Keywords},	
            pdfauthor={Authors},		
            bookmarksopen=false,
            pdfpagemode=UseNone]{hyperref}
            
\usepackage[margin = 1.2in]{geometry}
\usepackage[T1]{fontenc}
\usepackage[english]{babel}
\usepackage[latin1]{inputenc} 
\usepackage{mathtools}
\usepackage{amsmath}
\usepackage{amsfonts}
\usepackage{amsthm}
\usepackage{amssymb}
\usepackage{array}
\usepackage{algorithm} 
\usepackage{algorithmicx}
\usepackage{algpseudocode}
\usepackage{subcaption}
\usepackage{siunitx}

\usepackage{color}
\usepackage{url}
\usepackage{cleveref}
\usepackage{graphicx}
\usepackage{breakcites}
\usepackage{soul} 
\usepackage{enumitem}
\usepackage[title,titletoc,toc]{appendix}

\usepackage{makecell}

{\noindent {\textbf{Proof}:} }%
{\hfill $\Box$ \\[1ex] }

\newcommand{\bit}{\begin{itemize}}
	\newcommand{\eit}{\end{itemize}}
\newcommand{\ben}{\begin{enumerate}}
	\newcommand{\een}{\end{enumerate}}

\newcommand {\real} {\mathbb{R}}
\newcommand {\nat} {\mathbb{N}}

\DeclareMathOperator*{\Exp}{\mathbb{E}}
\DeclareMathOperator*{\Var}{\mathbb{V}\text{ar}}

\DeclareMathOperator*{\argmin}{arg\,min}%
\newcommand{\bA}{\ensuremath{\mathbf{A}}}

\newcommand{\bE}{\ensuremath{\mathbf{E}}}
\newcommand{\bF}{\ensuremath{\mathbf{F}}}

\newcommand{\bJ}{\ensuremath{\mathbf{J}}}

\newcommand{\bT}{\ensuremath{\mathbf{T}}}

\newcommand{\bW}{\ensuremath{\mathbf{W}}}
\newcommand{\bX}{\ensuremath{\mathbf{X}}}

\newcommand{\bZ}{\ensuremath{\mathbf{Z}}}

\newcommand{\ba}{\ensuremath{\mathbf{a}}}
\newcommand{\bb}{\ensuremath{\mathbf{b}}}
\newcommand{\bc}{\ensuremath{\mathbf{c}}}
\newcommand{\bd}{\ensuremath{\mathbf{d}}}

\newcommand{\bg}{\ensuremath{\mathbf{g}}}

\newcommand{\bj}{\ensuremath{\mathbf{j}}}

\newcommand{\br}{\ensuremath{\mathbf{r}}}

\newcommand{\bx}{\ensuremath{\mathbf{x}}}

\newcommand{\bz}{\ensuremath{\mathbf{z}}}

\newcommand{\cB}{\ensuremath{\mathcal{B}}}

\newcommand{\cD}{\ensuremath{\mathcal{D}}}

\newcommand{\cF}{\ensuremath{\mathcal{F}}}

\newcommand{\cU}{\ensuremath{\mathcal{U}}}

\usetikzlibrary{shapes.geometric, arrows, shadows, fit}

\usepackage{nomencl}
\makenomenclature

\renewcommand*{\nompreamble}{\begin{multicols}{2}}
\renewcommand*{\nompostamble}{\end{multicols}}
\begin{document}
	
	\begin{frontmatter}
		
		\title{Robust Design Optimization with Limited Data for Char Combustion}
				
		\author[ucsd]{{Yulin Guo}\corref{first}}
		
		\author[hyu]{{Dongjin Lee}\corref{first}}
		
		\author[ucsd]{Boris Kramer\corref{cor1}}
		\ead{bmkramer@ucsd.edu; +1 858-246-5327}
		\cortext[first]{These two authors contributed equally to this work.}
		\cortext[cor1]{Corresponding author}
						
		\address[ucsd]{Department of Mechanical and Aerospace Engineering, University of California San Diego, CA, United States}
  		\address[hyu] {Department of Automotive Engineering, Hanyang University, Seoul, South Korea}
		
		\begin{abstract}
            This work presents a robust design optimization approach for a char combustion process in a limited-data setting, where simulations of the fluid-solid coupled system are computationally expensive. We integrate a polynomial dimensional decomposition (PDD) surrogate model into the design optimization and induce computational efficiency in three key areas. First, we transform the input random variables to have fixed probability measures, which eliminates the need to recalculate the PDD's basis functions associated with these probability quantities. Second, using the limited  data available from a physics-based high-fidelity solver, we estimate the PDD coefficients via sparsity-promoting diffeomorphic modulation under observable response preserving homotopy regression. Third, we propose a single-pass surrogate model training that avoids the need to generate new training data and update the PDD coefficients during the derivative-free optimization. The results provide insights for optimizing process parameters to ensure consistently high energy production from char combustion.
		\end{abstract}	
		
		
		
\end{frontmatter}

\nomenclature[01]{$\mathbf{a}(t)$}{Potential sD-MORPH solution}
\nomenclature[02]{$\bA^+$}{Intermediate matrix for calculating $\breve{\bc}$, details in \cite{lee2024global}}
\nomenclature[03]{$\mathbb{A}^{N}$}{Subdomain of $\real^N$}
\nomenclature[04]{$\mathcal{B}^{N}$}{Borel $\sigma$-field on $\mathbb{A}^{N}$}
\nomenclature[05]{$\breve{\mathbf{c}}$}{sD-MOPRPH estimates of the PDD coefficients}
\nomenclature[06]{$\mathbf{c}_0$}{LASSO estimates of the PDD coefficients}
\nomenclature[07]{$c$}{Expansion coefficient}
\nomenclature[08]{$C_{\text{CO}}$, $C_{\text{O}_2}$, $C_{\text{H}_2\text{O}}$}{Mass concentrations of CO, O$_2$, H$_2$O, $\mathrm{kmole/m^3}$}
\nomenclature[09]{$C_p$}{Specific heat capacity of the gas mixture $\mathrm{J/(kg\cdot K)}$}
\nomenclature[10]{$C_{pi}$}{Specific heat capacity of a gas component, $\mathrm{J/(kg\cdot K)}$}
\nomenclature[11]{$\mathcal{D}$}{Design space}
\nomenclature[12]{$\mathbf{d}$}{Vector of design variables}
\nomenclature[13]{$d_p$}{Particle size, $\mathrm{m}$}
\nomenclature[14]{$D_o$}{Oxygen-nitrogen mixture diffusion coefficient}
\nomenclature[15]{$\bar{\mathbf{E}}$}{Intermediate matrix for calculating $\breve{\bc}$, details in \cite{lee2024global}}
\nomenclature[16]{$\mathcal{F}_{\bd}$}{$\sigma$-field on $\Omega_{\bd}$}
\nomenclature[17]{$\bar{\mathbf{F}}$}{Intermediate matrix for calculating $\breve{\bc}$, details in \cite{lee2024global}}
\nomenclature[18]{$\mathbf{g}$}{Vector of mean values}
\nomenclature[19]{$h(\cdot)$}{Thermal energy as a function of transformed random variables}
\nomenclature[20]{$h_{S,m}$}{$S$-variate, $m$-th order PDD approximation of $h$}
\nomenclature[21]{$\breve{\mathcal{K}}$}{Cost function in sD-MORPH regression}
\nomenclature[22]{$L$}{Number of PDD basis functions}
\nomenclature[23]{$m$}{Number of highest PDD order}
\nomenclature[24]{$m_{ci}$}{Unreacted char mass, $\mathrm{g}$}
\nomenclature[25]{$\dot{m}$}{Mass flow rate, $\mathrm{g/s}$}
\nomenclature[26]{$M_i$}{Mole fraction of a gas component, $\mathrm{g/mol}$}
\nomenclature[27]{$N$}{Number of random inputs}
\nomenclature[28]{$\Omega_{\bd}$}{Sample space}
\nomenclature[29]{$p_o$}{Oxygen partial pressure, $\mathrm{Pa}$}
\nomenclature[30]{$\mathbb{P}_{\bd}$}{A family of probability measures}
\nomenclature[31]{$P_s$}{Pressure constant}
\nomenclature[32]{$Q$}{Thermal energy, $\mathrm{J}$}
\nomenclature[33]{$R$}{Gas constant, $\mathrm{J/(mole\cdot K)}$}
\nomenclature[34]{$\real$}{Real numbers}
\nomenclature[35]{$\real_0^+$}{Non-negative real numbers}
\nomenclature[36]{$R_{\text{chem}}$}{Arrhenius kinetic rate, $\mathrm{s^{-1}}$}
\nomenclature[37]{$R_{\text{diff}}$}{Gas diffusion rate, $\mathrm{m^2/s}$}
\nomenclature[38]{$\mathbf{r}$}{Vector of deterministic variables used for transformation}
\nomenclature[39]{$S$}{Number of variate in PDD}
\nomenclature[40]{$Sh$}{Sherwood number}
\nomenclature[41]{$T$}{Temperature, $\mathrm{K}$}
\nomenclature[42]{$\bar{\mathbf{T}}$}{Intermediate matrix for calculating $\breve{\bc}$, details in \cite{lee2024global}}
\nomenclature[43]{$T_{\text{avg}}$}{Average temperature at the pressure outlet, $\mathrm{K}$}
\nomenclature[44]{$\mathbf{W}$}{Diagonal matrix in cost function $\breve{\mathcal{K}}$}
\nomenclature[45]{$w_1$, $w_2$}{Weighting factors in the objective function}
\nomenclature[46]{$\mathbf{X}$}{Uncertain input vector}
\nomenclature[47]{$y(\cdot)$}{Generic output function}
\nomenclature[48]{$\mathbf{Z}$}{Vector of transformed random variables}
\nomenclature[49]{$\mathbf{z}$}{Realization of $\bZ$}

\nomenclature[50]{$\beta$}{Empirical unitless exponent}
\nomenclature[51]{$\epsilon_{cp}$}{Maximum possible packing fraction}
\nomenclature[52]{$\epsilon_s$}{Empirical pressure constant}
\nomenclature[53]{$\lambda$}{Weight in cost function $\breve{\mathcal{K}}$}
\nomenclature[54]{$\mu_0$, $\sigma_0$}{Normalizing factors for the mean and standard deviation in the objective function}
\nomenclature[55]{$\mathbf{\Phi}$}{Intermediate matrix for calculating $\breve{\bc}$, details in \cite{lee2024global}}
\nomenclature[56]{$\Psi(\cdot)$}{Univariate orthonormal polynomial function}
\nomenclature[57]{$\tau$}{Collision stress, $\mathrm{Pa}$}

\section{Introduction} \label{sec:intro}
Biomass power plants convert organic materials, including wood, agricultural residues, and dedicated energy crops, into electricity and heat through combustion. The process includes several stages: drying (the initial phase that removes moisture from the biomass fuel), pyrolysis (the thermal decomposition of the biomass that produces volatile gases and solid char), and combustion (a chain of chemical reactions of the volatile gases and char that produces carbon dioxide and water vapor). The biomass combustion process is considered a sustainable energy generator, since the organic fuel is produced by absorbing CO\textsubscript{2} during plant growth. This process circulates CO\textsubscript{2} in the atmosphere instead of the carbon gas addition in the circulation~\cite{osman2023cost}. In contrast, the combustion of fossil fuels, which have been preserved over long periods, results in a net addition of CO\textsubscript{2} and is therefore deemed unsustainable, as it contributes to the global climate crisis. Improving combustion efficiency minimizes waste and gas emissions, making the process more sustainable overall. 

The complex and multiscale biomass combustion process is subject to uncertainties in fuel variability (type, moisture content, pre-treatment processes, etc.), operating conditions (temperature, oxygen availability, etc.), and emission measurement uncertainties~\cite{acp-19-8523-2019, steiner2024particulate}, which makes it difficult to control and predict the outcomes of biomass combustion. To achieve consistent energy production with biomass combustion, it is essential to optimize the process while taking into account these uncertainties. Reliability-based design optimization~\cite{tu1999new, zou2006direct, ren2016reliability} is one paradigm of design under uncertainty framework. By considering the probability of failure as a constraint, reliability-based design optimization aims to achieve high reliability of the optimal design. Another common design under uncertainty approach is robust design optimization (RDO)~\cite{youn2005performance, huang2006robust, beyer2007robust, huang2007analytical}, where a multi-objective optimization problem is formulated that minimizes a combination of the mean and variance of the performance function. It aims to achieve designs that perform consistently across different conditions. The RDO proposed herein allows one to maximize combustion efficiency and minimize fluctuations in the produced thermal energy. However, the optimization of the parameters related to the fuel pre-treatment and furnace operation requires repeated evaluations of the computationally expensive fluid-solid coupled combustion model. Surrogate models can replace the physics-based model to reduce the computational time for RDO, yet given that only a limited amount of high-fidelity training data can be afforded, we explore methods for training surrogate models in a limited-data setting.

The combustion process involves complex chemical reactions and heat transfer in a fluid-solid coupled model. The discrete element method (DEM) simulates the individual particle behavior through detailed interactions and dynamics in combustion \cite{xie2019coupling, xie2021study}. Another approach, the particle-in-cell (PIC) method \cite{snider2001incompressible}, models the gas using a Eulerian method while the solids are modeled using a Lagrangian method. Particles with identical physical properties are grouped to efficiently track their positions and trajectories. However, despite these model assumptions, the computational cost remains high. Conventional RDO evaluates statistical moments--such as mean and variance--in the objective function and constraints using quadrature rules~\cite{youn2005performance, huang2006robust}, which are only applicable when one has low-dimensional inputs. A straightforward  alternative is to evaluate statistical moments  by repeatedly sampling uncertain inputs according to their joint distribution and conducting physics-based simulations. Techniques such as importance sampling \cite{robert2013monte, doucet2000sequential} and Latin hypercube sampling \cite{helton2003latin, janssen2013monte} can reduce the number of samples required, but the necessity for expensive simulations remains. Surrogate models can be used within design optimization as cheaper-to-evaluate computational models, yet with the caveat of introducing model error. The surrogate models include, but are not limited to, polynomial chaos expansion~\cite{ren2013robust, shen2016polynomial} and polynomial dimensional decomposition (PDD) \cite{lee2023bi, lee2023multifidelity, rahman2008polynomial, rahman2009extended}, Kriging \cite{wang2021efficient, yang2021hybrid, jin2003use}, spline dimensional decomposition~\cite{LEE2022103218}, copula-based surrogates \cite{noh2009reliability}, reduced-order modeling~\cite{weickum2006multi}, and artificial neural networks~\cite{chatterjee2019critical}. To reliably optimize, the surrogate model needs to be accurate, which may require a substantial amount of high-quality training data from the costly physics-based model. Moreover, within the optimization, one may have to retrain the surrogate model as the design variables change. 

In this study, we formulate and solve a RDO problem for a char combustion process in a limited-data setting, since the multiscale fluid-solid combustion simulation is computationally expensive. Our contributions are as follows: (a) We formulate the optimization problem using transformed input random variables to further reduce the computational cost. The probability measures of the transformed input random variables are fixed during the design iterations, therefore, we eliminate the need to recalculate quantities associated with the probability quantities. (b) We adapt a sparsity-promoting diffeomorphic modulation under observable response-preserving homotopy (sD-MORPH) regression approach to accurately train a PDD surrogate model from limited data specifically for design optimization. We choose the PDD surrogate instead of other surrogates, such as polynomial chaos expansion, because we found PDD to better approximate the highly nonlinear quantity of interest. The PDD reduces the curse of dimensionality by effectively truncating basis functions that involve higher-order interactions among inputs. (c) We develop a single-pass surrogate model training process that estimates the objective function corresponding to the current design variable values with the PDD coefficients calculated in the initial design. There is no need to obtain input-output data from the expensive physics-based models at every optimization iteration. Although this is less accurate than updating the surrogate at every iteration, we demonstrate that the final optimized parameters achieve the desired significant increase in thermal energy from the combustion process. (d) Unlike previous works that relied on gradient-based optimization solvers, which is sensitive to noisy and nonsmooth outputs~\cite{lee2021robust,lee2023high}, the proposed RDO approach is more effective as the single-pass surrogate model training process leverages the advantages of a derivative-free design method. (e) We obtain Pareto solutions to the design variables, which provide insights into achieving a stable thermal production of the combustion process. In summary, we enhance the computational efficiency of the RDO problem by integrating a single-pass surrogate model training that involves training a PDD model with sD-MORPH regression. Additionally, we optimize the process by utilizing transformed input random variables. The RDO problem for char combustion is solved in a data-sparse setting which is common when working with computationally expensive high-fidelity solvers.

The rest of the paper is organized as follows. Section~\ref{section:problem_description} provides a high-level overview of the objective of the RDO problem, the relevant physics-based models, and the design optimization formulation. Section~\ref{section:char_combustion} introduces the char combustion model, the computational domain, the governing equations and the PIC method to numerically solve them, and defines the quantity of interest. Section~\ref{section:RDO} has all the details about the RDO problem, including the formulation of the design optimization problem and the surrogate model that works with limited data. Section~\ref{section:numerical_results} presents the numerical results of the RDO problem for the char combustion process. Section~\ref{section:conclusion} offers conclusions and an outlook towards future work.
%
\section{Problem description}\label{section:problem_description}
\begin{figure}[htbp]
    \centering 
        \begin{subfigure}{0.606\textwidth}{\includegraphics[width=\linewidth]{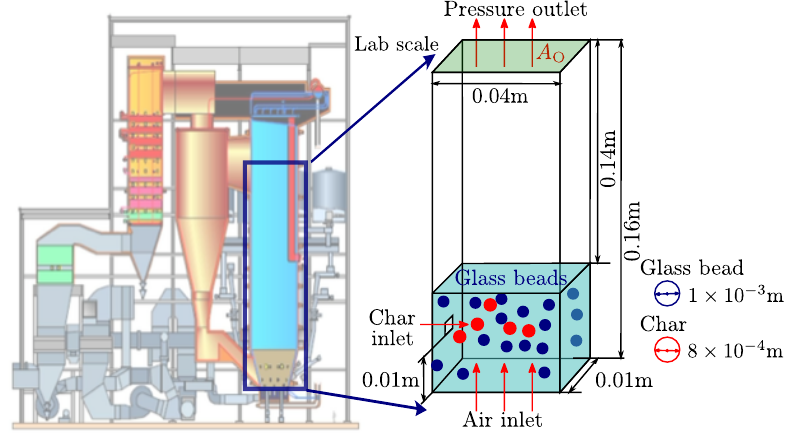}}\caption{\textcolor{black}{Geometry model}}\label{fig1:(a)}
        \end{subfigure} 
        \begin{subfigure}{0.386\textwidth}{\includegraphics[width=\linewidth]{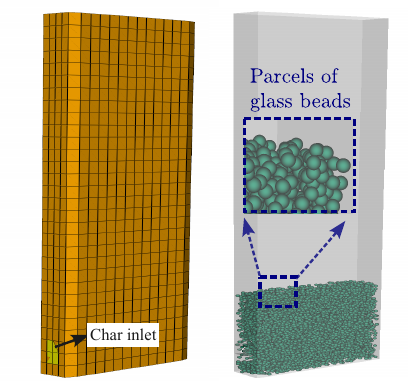}}\caption{Computational model}\label{fig1:(b)}
        \end{subfigure} 
        
        \caption{Fluidized bed furnace: (a) an industrial-scale model of the fluidized bed furnace (left) is scaled down to a lab-scale version (right),  which details the geometry and initial conditions with the diameters of glass beads and char particles, the air inlet, and char inlet; (b) in one simulation using the PIC method, we use 2,520 cells to discretize the fluid (left) and the parcel model uses 8,344 parcels to predict solid behavior (right).}\label{fig:1}
\end{figure}

Our goal is to optimize the char combustion process in a biomass power plant, specifically, to find the optimal operation parameters that maximize the produced thermal energy and minimize its variation in the presence of many uncertainties.
We simulate the char combustion process in a fluidized bed furnace model formulated on a lab-scale rectangular boiler (Figure~\ref{fig:1}) with the PIC method. The details of this computational model are discussed in Section~\ref{section:char_combustion}. The quantity of interest (QoI) of the char combustion process is the thermal energy $Q$ generated in a specified time interval, see~\eqref{equation:QoI} for a precise definition. It is a function of five uncertain inputs $\bX$ (Table~\ref{tab:random_inputs}), which include parameters related to the fuel pre-treatment and furnace operation, and depends on the design parameters.

The goal of the RDO problem is to maximize the thermal energy and minimize its variation, thus the objective function is chosen as a weighted sum of the first two statistical moments of $Q$. Specifically, we formulate the RDO problem as
    \begin{equation}
    \begin{aligned}
        \min_{\bd \in \cD \subseteq \real^{n}} & w_1\frac{\mu_0}{\mathbb{E}_{\bd}[Q(\bX;\bd)]} + w_2\frac{\sqrt{\mathbb{V}\text{ar}_{\bd}[Q(\bX; \bd)]}}{\sigma_0}, 
    \end{aligned}
    \label{eq:RDO_formulation}
    \end{equation}
where $\mu_0 \in \real \setminus \{0\}$ and $\sigma_0 \in \real \setminus \{0\}$ are two non-zero, real-valued scaling factors, which are chosen to better condition, such as normalizing, the objective function; $\bd$ is the vector of $n$ design variables; $\mathcal{D} \coloneqq \times^{n}_{k = 1}[d_{k,L}, d_{k,U}]$ is the design space, where $d_{k,L}$ and $d_{k,U}$ are the lower and upper bounds of the $k$th design variable. The design variables are the distribution parameters, in this case the mean values, of a subset of the random inputs $\bX$ that influence $Q$ the most. The subset is chosen based on variance-based global sensitivity analysis in our previous work~\cite{lee2024global}. Here, $\Exp_{\bd}[\cdot]$ and $\Var_{\bd}[\cdot]$ are the expectation and variance operators, respectively, with respect to the probability measure $f_{\bX}(\bx;\bd) {\rm d} \bx$, which depends on $\bd$. Also, $w_1 \in \real_0^+$ and $w_2 \in \real_0^+$ are two non-negative, real-valued weights such that $w_1 + w_2 = 1$. 

    
In order to efficiently calculate the moments, we formulate the optimization with transformed random variables (Section~\ref{section:RDO_formulation}) and use a PDD surrogate model (Section~\ref{section:PDD}) in place of the expensive physics-based simulations. Considering the computational cost of obtaining training samples, we train the PDD surrogate model using an sD-MORPH regression approach (Section~\ref{section:D-MORPH}). We expedite the optimization with a single-pass surrogate model training process (Section~\ref{section:single_pass}) that also leverages the advantages of a derivative-free solver.
%
\section{Computational model for char combustion}\label{section:char_combustion}
%
To predict the biomass combustion behavior, it is crucial to describe the behaviors of particles and gases in the combustor. To achieve this, one can use various methods such as the two-phase fluid method \cite{boileau2008investigation}, the DEM method \cite{xie2019coupling, xie2021study}, or the PIC method \cite{snider2001incompressible}. Considering both accuracy and efficiency, we select the PIC method to efficiently describe the behavior of millions of particles and gases. The computational domain and the PIC method are discussed in Section~\ref{section:Char_combustion_PIC_domain} and Section~\ref{section:Char_combustion_PIC_method}, respectively. We subsequently use the combustion predictions obtained from the PIC models to train a surrogate model for the QoI, which is defined in Section~\ref{section:QoI}. The simulation validation is discussed in Section~\ref{section:Char_combustion_validation}. 
\subsection{Computational domain} \label{section:Char_combustion_PIC_domain}
Figure~\ref{fig1:(a)} shows the furnace model that is a lab-scale rectangular boiler with specific dimensions (40 mm $\times$ 160 mm $\times$ 10 mm). The furnace contains stacked glass beads at the bottom to facilitate combustion. Glass beads are inert bed material, due to their excellent fluidization and heat transfer properties, high sphericity, and breaking resistance, they improve the mixing of solids, heat transfer, and reaction rate~\cite{KORKERD2024119262, cui2007fluidization} in char combustion. Char particles are fed into the boiler at a constant rate through the 2 mm $\times$ 2 mm char inlet on the left wall. The char particles react with oxygen in air flowing in from the inlet at the bottom at a constant rate. We consider five random inputs in this model, namely, the height of the stacked glass beads, the diameters of the glass beads and the char particles, and the inflow rates of the char and the air. Further details are discussed in Section~\ref{section:numerical_results}.

\subsection{Computational method} \label{section:Char_combustion_PIC_method}
%
We use a PIC-based numerical model in three spatial dimensions to predict the combustion behavior in the furnace. The PIC model treats the gas phase as a continuous medium using the Eulerian method. The model uses the Lagrangian method for the solids phase, grouping particles into parcels based on their physical properties (e.g., density, diameter). This simplification of the particle kinetics,  rather than using the full Newtonian equation, also increases computational efficiency. For more details on PIC, refer to \cite{lee2024global}. 

In PIC, a collisional stress model accounts for interactions among particles and interactions between particles and walls \cite{snider2001incompressible}. The collision stress is expressed as $\tau = P_s \epsilon_s^{\beta}/\max(\epsilon_{cp} - \epsilon_s,\ \alpha (1 - \epsilon_s))$, where $P_s = 1.0$ is the empirical pressure constant, $\epsilon_s$ is the particle volume fraction, $\beta = 2.0$ is the empirical unitless exponent, $\epsilon_{cp}$ is the maximum possible packing fraction for particles, and $\alpha = 10^{-9}$ $\ll 1$.

%
%
The combustion process generates two gas products, CO and CO\textsubscript{2}, from the heterogeneous (gas-solid) and homogeneous (gas-gas) reactions
    \begin{align*}
    \text{C(Solid)} + 0.5\text{O}_2\text{(Gas)} \rightarrow \text{CO(Gas)},\quad \quad \quad
    \text{CO(Gas)} + 0.5\text{O}_2\text{(Gas)} \rightarrow \text{CO}_2\text{(Gas)}.
    \end{align*}
The collisions among particles and between particles and the wall lead to the ash falling off the particles, following an Arrhenius kinetic rate and gas diffusion rate \cite{xie2021study}:
    \begin{align*}
    \frac{{\rm d}m_{ci}}{{\rm d}t} = -\pi d_p^2 p_0 \dfrac{1}{\frac{1}{R_{\text{diff}}} + \frac{1}{R_{\text{chem}}}},~ 
    R_{\text{diff}} = \frac{24 Sh D_o}{d_p R T_m},~ 
    R_{\text{chem}} = A_i \exp\left( -\frac{E_i}{R T_p} \right),~ 
    d_p = \left( \frac{6 m_p}{\pi \rho_p} \right)^{1/3}.
    \end{align*}
Here, $m_{ci}$ is the unreacted char mass [kg], $R_{\text{diff}}$ is the gas diffusion rate [m$^2$/s] and $R_{\text{chem}}$ is the Arrhenius kinetic rate [s$^{-1}$]. Also, $Sh$ is the Sherwood number [dimensionless], $R$ is the gas constant [J/(mol$\cdot$K)], and $T$ is the temperature [K]. Additionally, $p_o$ is the oxygen partial pressure [Pa] and $D_o$ is the oxygen-nitrogen mixture diffusion coefficient [m$^2$/s], while $d_p$ is the particle size [m] as it shrinks due to the mass loss.


The homogeneous reaction can be calculated by the law of mass action via the Arrhenius formula, proposed by Dryer and Glassman \cite{dryer1973high}, i.e.,
    \begin{align*}
    r_{\text{CO}} = 3.98 \times 10^{14} \exp\left( -\frac{1.67 \times 10^5}{R T_g} \right) C_{\text{CO}} C_{\text{O}_2}^{0.25} C_{\text{H}_2\text{O}}^{0.5},
    \end{align*}
where $C_{\text{CO}}$, $C_{\text{O}_2}$, and $C_{\text{H}_2\text{O}}$ are the mass concentrations [kg/m$^3$] of CO, O$_2$, and H$_2$O, respectively.


The computational mesh for the gas phase consists of 2,520 cells. Each cell has six degrees of freedom, which include three velocity components and three scalar variables (temperature, species concentrations, and pressure). The total number of degrees of freedom for the full model is 15,120. The parcel weight is set to 3 in this work, meaning one parcel contains three particles. As the height of the freeboard and diameters of the particles are random, the total number of particles varies in every simulation. In Figure~\ref{fig1:(b)}, there are 8,344 parcels representing the glass beads in the boiler model.
We use the open-source software MFiX (version 23.1.1) \cite{osti_1630426} with an MPI-based parallel computing solver on 15 CPUs (Intel Xeon W-3175X CPU \@ 3.10 GHz) for the combustion simulations.
\subsection{Quantity of interest} \label{section:QoI}
Our quantity of interest (QoI) is the total thermal energy generated in 10 seconds, computed as 
    \begin{equation}
    Q(\bX) = \int_{t=0}^{t=10} \dot{Q}(\bX, t) {\rm d}t = \int_{t=0}^{t=10} C_p(\bX, t) \times \dot{m}(\bX, t) \times T_{\text{avg}}(\bX, t) {\rm d}t,
    \label{equation:QoI}
    \end{equation}
where $T_{\text{avg}}(\bX, t)$ is the average temperature at the pressure outlet (at the top of the boiler model, c.f. Figure~\ref{fig:1}) at time $t$ and inputs $\bX$; $\dot{m}(\bX, t)$ is the mass flow rate at the outlet and $C_p(\bX, t)$ is the specific heat capacity of the mixture of O\textsubscript{2}, N\textsubscript{2}, CO, CO\textsubscript{2} and H\textsubscript{2}O, i.e., $C_p(\bX, t) = \sum_{i=1}^{5} C_{pi}(\bX, t) \times M_i(\bX, t)$. Here, $C_{pi}$ is the specific heat capacity of each gas component in the mixture and $M_i$ is the mole fraction of each gas component, such that $\sum_{i=1}^5 M_i=1$. The specific heat capacity $C_{pi}$ is a function of the gas constant, the molecular weight of the gas, and the temperature. Details are presented in Appendix~\ref{section:appendix_cp}. 

%
\subsection{Validation} \label{section:Char_combustion_validation}
We validate the PIC-based computation model with the DEM results from \cite[Fig. 2]{xie2021study}, focusing on the time evaluation of gas mass  fractions for CO and CO\textsubscript{2}. The root-mean-squared error between the DEM data and the data from our PIC-based solver for mass fractions at 5-second intervals over 40 seconds (see \cite{lee2024global}, Figure 6) is 2.73\%. This indicates that the PIC model is accurate and validated compared to the computational model used in \cite{xie2021study}.
%
\section{Robust design optimization problem and its solution with limited data} \label{section:RDO}
%
For optimal and consistent operation of the char combustion process in the fluidized bed furnace, we maximize the mean of the thermal energy while minimizing its variance to limit power fluctuations. The most influential design parameters were identified in our previous work through a global sensitivity analysis~\cite{lee2024global}; they are the mean values of the air inflow rate and the diameter of the glass bead particle. However, estimating the mean and variance of the thermal energy requires evaluating hundreds to thousands of samples via costly combustion simulations. While a surrogate model can replace the high-fidelity simulation, a typical optimization requires the surrogate model to be trained repeatedly. To circumvent this computational burden, we propose a single-pass surrogate model training process as well as two other efficiency boosts.
  
Section~\ref{section:formulation} defines the notations and the char combustion problem setup. Section~\ref{section:RDO_formulation} discusses RDO formulations with transformed random variables. In Section~\ref{section:PDD} we briefly review the PDD surrogate model from \cite{rahman2008polynomial, rahman2009extended} to keep the paper self consistent.  Sections~\ref{section:moments}~and~\ref{section:D-MORPH} discuss the calculation of the PDD coefficients through sD-MORPH regression with limited data, and using the coefficients to calculate the statistical moments of the QoI. Finally, the single-pass training of PDD and the complete algorithm for RDO are discussed in Sections~\ref{section:single_pass}~and~\ref{section:complete_algorithm}.
\subsection{Problem setup and definitions} \label{section:formulation}
We define random variables within a measurable space, consisting of a sample space and associated probability measures. As the measurable space depends on the design vector $\mathbf{d}$, we include $\mathbf{d}$ in the notation for the associated variables and operators hereafter. We denote the measurable space by $(\Omega_{\bd},\ \cF_{\bd})$, where $\Omega_{\bd}$ represents a sample space and $\cF_{\bd}$ is a $\sigma$-field on $\Omega_{\bd}$. A family of probability measures $\{\mathbb{P}_{\bd}: \ \cF_{\bd} \rightarrow [0, 1]\}$ is defined over $(\Omega_{\bd},\ \cF_{\bd})$. 

The input random vector $\bX \coloneqq (X_1,\ X_2,\ \ldots,\ X_N)^{\top}$ describes uncertainties in the system, $N \in \nat$. It is an $\mathbb{A}^{N}$-valued input random vector with $\cB^{N}$ representing the Borel $\sigma$-field on $\mathbb{A}^{N} \subseteq \real^N$, $(\Omega_{\bd},\ \cF_{\bd}) \rightarrow (\mathbb{A}^{N},\ \cB^{N})$. The probability law of $\bX$ is completely defined by a family of the joint probability density functions (PDFs) $\{f_{\bX}(\bx;\bd): \bd \in \real^{N}, \bd \in \cD \}$ that are associated with probability measures $\{\mathbb{P}_{\bd}:\ \bd \in \cD\}$. 


In this study, we define the random variables as the height of the freeboard ($X_1$), the air inflow rate ($X_2$), the diameter of the glass bead particle ($X_3$), the diameter of the char particle ($X_4$), and the char mass inflow rate ($X_5$) so $N = 5$. The $n$-dimensional design vector is defined as $\bd \coloneqq (d_1,\ d_2,\ \ldots, \ d_n)^{\top} \in \cD$, where $\cD$ is the design space. While the design variable $d_k$ could represent any distribution parameter such as the standard deviation, in this work, we assume that the design variable is the mean of the random variable, i.e., $d_k=\mathbb{E}[X_{i_k}]$, $k=1,\ldots,n$ and $i_k\in\{1,\ldots,N\}$. Note that only the mean values of the most influential random variables are chosen as the design variables. Based on variance-based global sensitivity analysis in~\cite{lee2024global}, they are the mean values of the air inflow rate and the diameter of the glass bead particle, so $n=2$. We assume $Q(\bX;\bd)$ is a real-valued, square-integrable, measurable transformation on $(\Omega_{\bd},\ \cF_{\bd})$, describing the quantity of interest.
\subsection{Robust design optimization with transformed random variables} \label{section:RDO_formulation}
%
In this work, the RDO aims to maximize the mean and minimize the standard deviation of the thermal energy output $Q(\bX;\bd)$ simultaneously. We achieve this by minimizing the bi-objective function as described in~\eqref{eq:RDO_formulation}. 

%
Throughout the design iterations, we transform the original random input $\bX$ to $\bZ$ such that each transformed $Z_i,\ i = 1, 2, \ldots, N,$ has a fixed mean value. The RDO problem formulation can be easily refactored using $\bZ$. The transformed $\bZ$ is used for the PDD basis, which eliminate the need to recalculate probability measure-associated quantities. This is particularly beneficial for estimating the mean and the standard deviation of the QoI via the PDD surrogate, as discussed in the Section~\ref{section:PDD}.
We define $\bZ = (Z_1, Z_2, \ldots, Z_N)^{\top}$ as an $N$-dimensional vector of new random variables obtained by scaling $\bX$ as
    \begin{equation}
        \bZ \coloneqq \text{diag}[r_1, r_2, \ldots, r_N]\bX,
        \label{eq:Z_from_X}
    \end{equation}
where $\br \coloneqq (r_1, r_2, \ldots, r_N)^{\top}$ is an $N$-dimensional vector of deterministic variables. Corresponding to the $n$-dimensional design vector, we denote $(Z_{i_1}, Z_{i_2}, \ldots, Z_{i_n})^{\top}$ as a subvector of $\bZ$, where the $i_k$th random variable $Z_{i_k}$ corresponds to the $i_k$th original random variable $X_{i_k}$. The mean of $Z_{i_k}$ is $\Exp_{\bd}[Z_{i_k}] = d_k r_{i_k} =: g_k, k = 1, 2, \ldots, n$. The PDF of $\bZ$ is 
    \[
    f_{\bZ}(\bz; \bg) = |\bJ| f_{\bX}(\bx; \bd) = \left| \frac{1}{r_1 r_2 \cdots r_N} \right| f_{\bX}(\bx; \bd) = \left| \frac{1}{r_1 r_2 \cdots r_N} \right| f_{\bX}\left(\text{diag}\left[\frac{1}{r_1}, \frac{1}{r_2}, \ldots, \frac{1}{r_N}\right]\bz; \bd\right)
    \]
supported on $\bar{\mathbb{A}}^{N} \subseteq \real^N$. The absolute value of the determinant of the Jacobian is $\left|\bJ\right| = \left|\partial \bx/\partial \bz \right| = 1/(r_1 r_2 \ldots r_N)$ and $\bg := (g_1,g_2,\ldots,g_n)^{\top} \in \real^n$ is an $n$-dimensional vector of mean values. We choose $g_k$ to be \textit{one}. Thus, $\bg = \textbf{1}$ and at every iteration, $r_{k} = 1/d_{k}$.


We define $h(\bZ; \br) \coloneqq Q(\bX)$, so we can reformulate \eqref{eq:RDO_formulation} to
    \begin{equation}
    \begin{aligned}
        \min_{\bd \in \cD \subseteq \real^{n}} & w_1\dfrac{\mu_0}{\mathbb{E}_{\bg(\bd)}[h_0(\bZ;\br)]} + w_2 \dfrac{\sqrt{\mathbb{V}\text{ar}_{\bg(\bd)}[h_0(\bZ;\br)]}}{\sigma_0}
    \end{aligned}
    \label{eq:Tranformed_RDO_formulation}
    \end{equation}
where $\Exp_{\bg(\bd)}$ and $\Var_{\bg(\bd)}$ are the expectation and variance operators, respectively, with respect to the probability measure $f_{\bZ}(\bz;\bg) {\rm d} \bz$, which depends on $\bd$. For brevity, the subscript $\bg(\bd)$ will be denoted by $\bg$ in the rest of the paper.
The transformed formulation in \eqref{eq:Tranformed_RDO_formulation} is expressed in terms of the transformed input random variables $\bZ$ and is equivalent to the original formulation in \eqref{eq:RDO_formulation}. 
\subsection{Polynomial dimensional decomposition (PDD) surrogate model} \label{section:PDD}
Any square-integrable function $y(\cdot)$ defined on the probability space $(\Omega, \cF, \mathbb{P})$ can be represented with an infinite Fourier-polynomial expansion known as polynomial dimensional decomposition (PDD)~\cite{rahman2008polynomial, rahman2009extended}, given by
    \begin{equation}
    y(\bX) = y_0 + \sum_{\emptyset \neq \cU \subseteq \{1, \ldots, N\}} \sum_{\bj_{\cU} \in \nat^{|\cU|}} c_{\cU, \bj_{\cU}} (\bx) \Psi_{\cU, \bj_{\cU}}(\bX_{\cU})
    \end{equation}
    \begin{equation*}
    c_{\cU, \bj_{\cU}}(\bx):=\int_{\mathbb{A}_N} y(\bx) \Psi_{\cU, \bj_{\cU}} (\bx_{\cU}) f_{\bX} (\bx) d\bx, \quad
    \Psi_{\cU, \bj_{\cU}} (\bX_{\cU}) := \prod_{i \in \cU} \Psi_{i, j_i} (X_i),
    \end{equation*}
where $c_{\cU, \bj_{\cU}}(\bx)$ is the expansion coefficient, $\Psi_{\cU, \bj_{\cU}} (\bX_{\cU})$ is the multivariate orthonormal polynomial, $\Psi_{i,j_i}$ is a univariate orthonormal polynomial in $X_i$ of degree $j_i$ and it is consistent with the probability measure $f_{X_i}(x_i)dx_i$. The full PDD has infinite dimensions. In practice, we truncate the infinite-dimensional PDD by retaining only the degrees of interaction among input variables that are less than or equal to $S$, which can be determined by the number of sensitive inputs. As a result, the polynomial expansions with degree $m$, $S\leq m < \infty$ are preserved. The $S$-variate, $m$th-order PDD approximation to $h(\bZ; \br)$, $h_{S, m} (\bZ;\br)$ can be expressed as
    \begin{equation}
    h_{S, m} (\bZ;\br) := h_0(\br) + \sum_{\substack{\emptyset \neq \cU \subseteq \{1, \ldots, N\}\\{1 \leq |\cU| \leq S}}} \sum_{\substack{\bj_{\cU} \in \mathbb{N}^{|\cU|} \\ {|\cU| \leq |\bj_{\cU}| \leq m}}} c_{\cU, \bj_{\cU}}(\br) \Psi_{\cU, \bj_{\cU}} (\bZ_{\cU};\bg) \approxeq h(\bZ;\br).
    \label{eq:PDD}
    \end{equation}
Here, $\Psi_{i,j_i}$ is a univariate orthonormal polynomial in $Z_i$ of degree $j_i$, which is consistent with the probability measure $f_{Z_i}(z_i;\bg){\rm d} z_i$.


In \eqref{eq:PDD}, we can arrange the elements of the basis in any order such that $\{\Psi_{\cU, \bj_{\cU}} (\bZ_{\cU};\bg) : 1 \leq |\cU| \leq S, \, |\cU| \leq \bj_{\cU} \leq m\} = \{\Psi_2 (\bZ;\bg), \ldots, \Psi_L (\bZ;\bg)\}, \, \Psi_1 (\bZ;\bg) = 1$, where $\Psi_i(\bZ;\bg)$ is the $i$th basis function in the truncated PDD approximation and $L_{N,S,m} = 1 + \sum_{s = 1}^S \binom{N}{s} \binom{m}{s}.$
We can rewrite the PDD surrogate in terms of $\Psi_i(\bZ;\bg)$ as
    \begin{equation}
    h_{S, m} (\bZ;\br) = \sum_{i = 1}^{L_{N,S,m}}  c_i(\br) \Psi_i (\bZ;\bg).
    \label{eq:PDD_approximation}
    \end{equation}
\subsection{Statistical moment analysis via the PDD surrogate}\label{section:moments}
It is advantageous to use PDD coefficients for estimating statistical moments, such as the mean and the variance, because these moments are explicitly expressed in terms of coefficients due to the orthogonality of the basis functions~\cite{sudret2008global, lee2024global}. We use the $S$-variate, $m$th-order PDD surrogate $h_{S,m}(\bZ;\br)$ to estimate the mean and variance of $h(\bZ;\br)$. The orthonomality of the basis functions allows for the analytical formulation of the mean of $h(\bZ,\br)$
    \begin{equation}
    \mathbb{E}_{\bg}[h_{S,m}(\bZ;\br)]=c_1(\br)=\mathbb{E}_{\bg}[h(\bZ;\br)]\label{mean}
    \end{equation}
    and the variance of $h(\bZ,\br)$
    \begin{align}
    \mathbb{V}\text{ar}_{\bg}[h_{S,m}(\bZ;\br)]=\sum_{i=2}^{L}c_i^2(\br) \leq \mathbb{V}\text{ar}_{\bg}[h(\bZ;\br)].\label{variance}
    \end{align}
This requires an accurate computation of the expansion coefficients to obtain good estimates of the mean and variance of $h(\bZ;\br)$. To achieve this, particularly when the training data set is limited, we use a regression method called sD-MORPH regression~\cite{lee2024global}, which is detailed in the following section. 
\subsection{Regression problem and solution from limited data} \label{section:D-MORPH}
The coefficients $\bc=(c_1,c_2,\ldots,c_L)^{\intercal}\in\mathbb{R}^L$ of the PDD surrogate in \eqref{eq:PDD_approximation} can be obtained by solving a linear system 
    \begin{align}
    \underbrace{
    \begin{bmatrix}
    \Psi_1(\bz^{(1)};\bg) & \cdots &  \Psi_{L}(\bz^{(1)};\bg) \\
    \vdots                   & \ddots &  \vdots                           \\
    \Psi_1(\bz^{(M)};\bg) & \cdots &  \Psi_{L}(\bz^{(M)};\bg)
    \end{bmatrix}}_{=:\bA } \underbrace{\begin{bmatrix} c_1 \\ c_2 \\ \vdots \\ c_L  \end{bmatrix}}_{\bc}
    = \underbrace{\begin{bmatrix} h(\bz^{(1)};\br) \\ h(\bz^{(2)};\br) \\\vdots \\ h(\bz^{(M)};\br) \end{bmatrix}}_{=:\bb}.
    \label{linear}
    \end{align}
To solve this with limited data ($M<L$), it is common to use LASSO regression, which penalizes the $l^1$ norm of the regression coefficients, i.e., 
    \begin{align}
    \bc_0 = \argmin_{\bc\in\mathbb{R}^L}\left\{ (\bb-\bA\bc)^{\intercal}(\bb - \bA\bc) + k \sum_{i=1}^L|c_i|\right \}, 
    \label{lasso}
    \end{align}
where $k$ is a positive real number. The second term of \eqref{lasso} is a regularization term that penalizes the $l^1$ norm of the PDD's expansion coefficients, producing sparse solutions for underdetermined systems. However, LASSO regression introduces a bias to reduce variance and cannot infer more non-zero coefficients than the number of training samples.

In \cite{lee2024global}, we proposed a sparsity-promoting D-MORPH regression to train the PDD surrogate model in limited data. This regression leverages the sparsity from LASSO regression and iteratively improves its accuracy by combining it with the D-MORPH solution. While the previous work focuses on global sensitivity analysis to reduce the number of design variables, this work presents a RDO method via the PDD surrogate trained with sD-MORPH that uses the same training data in~\cite{lee2024global}. In sD-MORPH, we define a cost function that measures the difference between a potential sD-MORPH solution $\ba(t)\in \mathbb{R}^L$ and the LASSO estimates $\bc_0\in\mathbb{R}^L$ to obtain an initial sD-MORPH solution 
    \begin{align}\label{prior_costf}
    \bc_1=\argmin_{t\in\mathbb{R}}\left\{\dfrac{1}{2}(\ba(t)-\bc_0)^{\intercal}(\ba(t)-\bc_0)\right\}.
    \end{align}

Next, we augment the cost function by including the $l^2$ norm between the potential sD-MORPH solution and the solution obtained from the previous iteration. This sD-MORPH regression then minimizes the augmented cost function $\breve{\mathcal{K}}(\ba(t))$, i.e., 
    \begin{align}\label{new_costf}
    \breve{\bc}^{(i)} = \argmin_{t\in\mathbb{R}}\left\{\breve{\mathcal{K}}(\ba(t))=\dfrac{\lambda}{2}\left(\ba(t)-\bc_0\right)^{\intercal}\bW\left(\ba(t)-\bc_0\right)+\dfrac{1-\lambda}{2}\left(\ba(t)-\bc^{(i-1)}\right)^{\intercal}\bW\left(\ba(t)-\bc^{(i-1)}\right)\right\}
    \end{align}
with a non-negative real-valued weight $\lambda\in[0,1]$. As a result of the sD-MORPH process (see Section 3.2.2 in \cite{lee2024global}), the final sD-MORPH solution becomes 
    \begin{equation}\label{sdmorph} 
    \begin{split}
    \breve{\bc}^{(i)}=&\bar{\bF}^{(i-1)}_{L-r}\left(\bar{\bE}^{(i-1)\intercal}_{L-r}\bar{\bF}^{(i-1)}_{L-r}\right)^{-1}\bar{\bE}^{(i-1)\intercal}_{L-r}\bA^+\bb\ +
    \\&\bar{\bF}^{(i-1)}_{r}\left(\bar{\bE}^{(i-1)\intercal}_{r}\bar{\bF}^{(i-1)}_{r}\right)^{-1}\bar{\bE}^{(i-1)\intercal}_{r}\left(\bar{\bT}_r^{(i-1)}\right)^{-1}\mathbf{\Phi}\mathrm{diag}\left({\breve{\bc}^{(i-1)}}\right)^{-1}\left(\lambda\bc_0+(1-\lambda){\bc_1}^{(i-1)}\right),
    \end{split}
    \end{equation}
where $\bar{\mathbf{E}}^{(i)}_r$ and $\bar{\mathbf{E}}^{(i)}_{L-r}$, $\bar{\mathbf{F}}^{(i)}_r$ and $\bar{\mathbf{F}}^{(i)}_{L-r}$ are constructed from the first $r$ and the last $L-r$ columns of matrices $\bar{\mathbf{E}}^{(i)}$ and $\bar{\mathbf{F}}^{(i)}$, respectively, generated from the singular value decomposition
    \begin{equation}
    \mathbf{\Phi}\bW^{(i-1)}=\bar{\mathbf{E}}^{(i-1)}
    \begin{bmatrix}
    \bar{\mathbf{T}}_r^{(i-1)} & \mathbf{0}\\
    \mathbf{0} & \mathbf{0}
    \end{bmatrix}
    \bar{\mathbf{F}}^{(i-1)\intercal},
    \end{equation}
with $\bar{\mathbf{T}}^{(i-1)}$ representing an $r\times r$ diagonal matrix including nonzero singular values. Here, $\bW^{(i-1)}=\text{diag}[0,1/(\breve{c}_2^{(i-1)}+\epsilon),\ldots,1/(\breve{c}_{L}^{(i-1)}+\epsilon)]$, where its first element is \emph{zero} and the remaining elements are the reciprocal of  $\breve{c}_j^{(i-1)}$ for $j=2,\ldots,L$ with $\epsilon \ll 1$. 
%


\subsection{Single-pass surrogate model training process} \label{section:single_pass}
In the RDO process, the PDD surrogate model should ideally be recomputed using a new input-output dataset whenever the design variables change. Although the sD-MORPH regression reduces the required number of input-output data, the cumulative number of simulations can become prohibitive, especially when the number of iterations is large. This issue is especially pronounced when using a derivative-free optimization solver. 

We use a single-pass surrogate model training process to reduce the total number of high-fidelity solutions needed. In contrast to~\cite{lee2021robust}, which uses a gradient-based sensitivity method, our single-pass surrogate model training uses a non-gradient-based method that better handles the nonlinear objective function. To use the single-pass training process, we assume that (a) an $m$th-order PDD approximation $h_{S,m}(\bZ; \br)$, see~\eqref{eq:PDD_approximation}, at the initial design is valid for all possible designs; (b) the PDD coefficients for a new design, computed by recycling PDD generated for the previous design, remain accurate. 

We denote the transformation vectors associated with the previous and current designs with $\br$ and $\br'$, respectively. To satisfy the two assumptions, we use a higher order of $m=11$ than $m=9$ in our previous work~\cite{lee2024global} for the PDD and limit the design space to be represented by the PDD model computed at the initial design. We first compute the PDD coefficients $c_i(\br)$, $i = 1, 2, \ldots, L_{S,N,m}$, collectively denoted as $\bc(\br)$, for the initial design using the input-output data $\{\bz^{(l)}, h(\bz^{(l)}; \br)\}_{l = 1}^{L_{S,N,m}}$. For the next design in the optimization process, we modify the input data $\{\bz^{(l)}\}_{l = 1}^{L_{S,N,m}}$ to $\{\bz'^{(l)}\}_{l = 1}^{L_{S,N,m}}$, where $\bz'^{(l)} = \text{diag}\left[ \frac{r_1}{r_1'}, \frac{r_2}{r_2'}, \ldots,  \frac{r_N}{r_N'}\right]\bz^{(l)}$. We then represent the QoI at the updated design $\br'$ in terms of $\bz'$ as 

\begin{equation}
    \begin{aligned}
        h\left( \bz^{(l)} ; \br' \right) 
        :&= h \left( \text{diag} \left[ \frac{1}{r'_1}, \frac{1}{r'_2}, \ldots, \frac{1}{r'_N} \right] \bz^{(l)} \right)
        \\&= h \left( \text{diag} \left[ \frac{1}{r_1}, \frac{1}{r_2},\ldots, \frac{1}{r_N} \right] \text{diag} \left[ \frac{r_1}{r'_1}, \frac{r_2}{r'_2}, \ldots, \frac{r_N}{r'_N} \right] \bz^{(l)} \right)
        \\&= h \left( \text{diag} \left[ \frac{1}{r_1}, \frac{1}{r_1}, \ldots, \frac{1}{r_N} \right] \bz'^{(l)} \right) 
        =: h \left( \bz'^{(l)} ; \br \right).
    \end{aligned}
\end{equation}

This modification allows us to approximate the output function $h(\bz^{(l)};\br')$ at the updated design $\br'$ via PDD computed at the initial design $\br$ as
    \begin{equation}
        h(\bz^{(l)}; \br') = h(\bz'^{(l)}; \br) \approx \sum_{i=1}^{L_{S,N,m}} c_i(\br) \Psi_i(\bz'^{(l)}; \bg)=:\tilde{h}(\bz^{(l)};\br'),
        \label{eq:h_approximation_by_prev_design}
    \end{equation}
which does not require updating the PDD, only evaluating the PDD at the new input $\bz'$. As a result, at every iteration, we estimate the mean and variance with $\bc(\br')$, resulting in little extra computational cost.

\begin{figure}[ht]
\centering 
{\begin{tikzpicture}[auto,
    block_center/.style ={rectangle, rounded corners, draw=black, thick,
      text width=18em, text centered, inner sep=6pt},
    block_left/.style ={rectangle, rounded corners, draw=black, thick,
      text width=16em, text centered,  text width = 6em, inner sep=6pt},
    line/.style ={draw, thick, -latex', shorten >=0pt}][h]
    
    \matrix [column sep = 1cm, row sep = 0.5cm]{
    & \node (r1) [block_center] {Approximate output $h(\bz^{(l)};\br')$ \\
    using $\bc(\br),$ see~\eqref{eq:h_approximation_by_prev_design}};\\
    \node (l1) [block_left] {Optimizer}; & 
    \node (r2) [block_center] {Construct the input-output data set \\ $\{\bz^{(l)}, h(\bz^{(l)};\br'\}^{M}_{l=1},\ M<L_{S,N,m}$};\\
    & \node (r3) [block_center] {Compute $\bc(\br')$, the PDD \\ coefficients at $\br'$, using sD-MORPH};\\};

    \coordinate (invisible) at ($(l1) - (6em, 0)$);
    
    \begin{scope}[every path/.style = line]
        \path (l1) |-  (r1) node[above, near end]{$\br'$};
        \path (r3) -| (l1);
        \path (r1) -- (r2);
        \path (r2) -- (r3);
        \path (invisible) --node {$\br$} (l1);
    \end{scope}
\end{tikzpicture}}
\caption{A flow chart for a single-pass surrogate model training process that repeatedly computes PDD for the updated design $\br'$ using the PDD computed at $\br$ with a sD-MORPH.}\label{fig:flowchart_single_pass}
\end{figure}
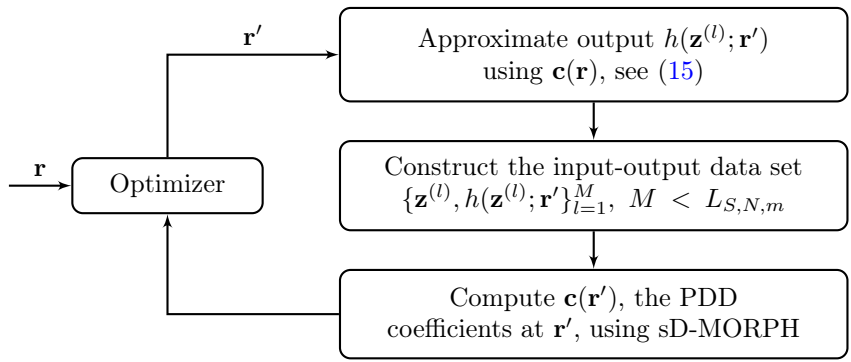

We select the Nelder-Mead method~\cite{mckinnon1998convergence} as the derivative-free optimizer. Combustion processes often involve complex chemical reactions, resulting in objective functions that are highly non-linear, non-smooth, and possibly discontinuous. The Nelder-Mead algorithm does not require gradient information and can handle such output complexities more effectively. Figure~\ref{fig:flowchart_single_pass} shows a flowchart for the proposed single-pass training process that includes training the PDD surrogate by the sD-MORPH regression during RDO. 

\subsection{Complete algorithm for robust design optimization}\label{section:complete_algorithm}

\begin{algorithm}
\caption{RDO with limited data}
\begin{algorithmic}[1]
\State \textbf{Initialize:} Set $\bd=\bd_0$, $i=1$
\State Sample $\bX \sim f_{\bX}(\bx;\bd)\mathrm{d}\bx$
\State Compute $\{\bx^{(l)}, Q(\bx^{(l)};\bd)\}_{l=1}^M$, where $M < L_{S,N,m}$

\State Transform $\bX$ to $\bZ$ according to~\eqref{eq:Z_from_X}
\State Construct $\{\bz^{(l)}, h(\bz^{(l)};\br)\}^{M}_{l=1}$

\State Compute $\bc(\br)$, coefficients of $h_{S,m}(\bz;\br)$

\State Estimate the mean using~\eqref{mean} and variance using~\eqref{variance} of the output $h(\bz;\br)$

\State Derivative-free optimization for $\bd_{*}^{\{1\}}$

\While{not converged}

    \State  Approximate output $h(\bz^{(l)};\br')$ using PDD composed of $\bc(\br)$ using~\eqref{eq:h_approximation_by_prev_design}
    
    \State Construct $\{\bz^{(l)}, \tilde{h}(\bz^{(l)};\br')\}^{M}_{l=1}$
    
    \State Compute $\bc(\br')$, coefficients of $h_{S,m}(\bz;\br')$

    \State Estimate the mean using~\eqref{mean} and variance using~\eqref{variance} of the output $h(\bz;\br')$
    
    \State Derivative-free optimization for $\bd_{*}^{\{i\}}$
    
    \State $i \gets i + 1$
    
\EndWhile

\State \textbf{Return:} Set solution $\mathbf{d}^* = \mathbf{d}^{\{i\}}$
\end{algorithmic}\label{alg:RDO}

\end{algorithm}

\begin{figure}[ht]
    \centering
    \tikzstyle{arrow} = [thick,->,>=stealth]

\begin{tikzpicture}[auto,
    block_1/.style = {rectangle, rounded corners, draw=black, thick,
      text width=26em, text centered, inner sep=6pt},
    block_2/.style = {rectangle, draw=black, thick,
      text width=26em, text centered, inner sep=6pt},
    block_3/.style = {rectangle, draw=black, thick,
      text width=10em, text centered},
    block_4/.style = {rectangle, draw=black, thick,
      text width=17em, text centered, inner sep=6pt},
    block_5/.style = {rectangle, draw=black, thick,
      text width=22em, text centered, inner sep=6pt},
    decision/.style = {diamond, draw=black, thick},
    cloud/.style = {draw, ellipse, node distance=3cm, minimum height=2em},
    line/.style = {draw, thick, -stealth', shorten >=0pt},
    back group/.style = {fill = blue!10, rounded corners, draw = black!50, dashed, inner xsep = 15pt, inner ysep=10pt}]

    \centering
    \node (l1) [block_1] {Initialization: set $\bd=\bd_0$, $i=1$};
    \node (l2) [block_2, anchor = north] at ($(l1.south) + (0, -1.5em)$) {
    Sample $\bX$ $\sim$ $f_{\bX}(\bx;\bd)\mathrm{d}{\bx}$, compute physics-based model input-output data $\{\bx^{(l)}, Q(\bx^{(l)};\bd)\}_{l=1}^{M},\quad M<L_{S,N,m}$};
    \node (l3) [block_2, anchor = north] at ($(l2.south) + (0, -1.5em)$) {
    Transform $\bX$ to $\bZ$, see~\eqref{eq:Z_from_X}, construct $\{\bz^{(l)}, h(\bz^{(l)};\br)\}^{M}_{l=1}$};
    \node (l4) [block_3, anchor = north] at ($(l3.south) + (0em, -4em)$) {Compute $\bc(\br')$, coefficients of $h_{S,m}(\bz;\br')$};
    \node (r4) [block_4, anchor = south west] at ($(l4.south east) + (1.5em, 0em)$) {
    Approximate output $h(\bz^{(l)};\br')$ using PDD composed of $\bc(\br)$, see~\eqref{eq:h_approximation_by_prev_design}, construct $\{\bz^{(l)}, \tilde{h}(\bz^{(l)};\br')\}^{M}_{l=1}$};
    \node (ll4) [block_3, anchor = east] at ($(l4.west) + (-1.5em, 0em)$) {Compute $\bc(\br)$, coefficients of $h_{S,m}(\bz;\br)$};
    \node (l5) [block_5, anchor = north west] at ($(ll4.south west) + (0, -1.5em)$) {Estimate the mean, see~\eqref{mean} and variance, see~\eqref{variance} of the output $h(\bz;\br)$ or $h(\bz;\br')$}; 
    \node (l6) [block_5, anchor = north west] at ($(l5.south west) + (0, -1.5em)$) {Derivative-free optimization for $\bd_{*}^{\{i\}}$, $i=i+1$};
    \node (r6) [decision, anchor = north] at ($(r4.south) + (0em, -2em)$) {\small Converge?};
    \node (r7) [cloud, anchor = north] at ($(r6.south) + (0, -3em)$) {
    Set solution $\bd^* = \bd_{*}^{\{i\}}$};

    \draw[arrow] (l1.south) -- (l2.north);
    \draw[arrow] (l2.south) -- (l3.north);
    \draw[arrow] ($(ll4.north) + (0, 4em)$) -- (ll4.north);
    \draw[arrow] (l4.south) --++ (0, -1.5em);
    \draw[arrow] (ll4.south) --++ (0, -1.5em);
    \draw[arrow] (l5.south) -- (l6.north);
    \draw[arrow] (r4.west) -- ($(r4.west) + (-1.5em, 0)$);

    \draw[thick] (l6.east) -- ($(l6.east) + (3.5em, 0em)$);
    \draw[arrow] ($(l6.east) + (3.5em, 0em)$) |- (r6.west);

    \draw[arrow] (r6.south) -- node [left, near end] {Yes} (r7.north);
    \draw[arrow] (r6.north) -- node [right, midway] {No} (r4.south);

    \begin{scope}[on background layer]
        \node (bk) [back group] [fit = (l4) (ll4) (l5) (l6) (r4) (r6)] {};
    \end{scope}

    \node[above, yshift = 1.5em, xshift = -9em, text = blue, right] at (l4.north) {Single-pass training process};
    
\end{tikzpicture}
    \caption{Flowchart for the RDO with limited data.}
    \label{fig:flowchart_algorithm}
\end{figure}

Algorithm~\ref{alg:RDO} presents the pseudocode of the proposed RDO method in the limited data setting. Figure~\ref{fig:flowchart_algorithm} illustrates the corresponding flowchart. This algorithm begins with initialization that includes setting an initial design vector $\mathbf{d}=\mathbf{d}_0$ and all parameters (e.g., $N$, $S$, $m$ used in PDD and termination criteria $|\epsilon| \ll 1$ used in design process). We then generate input samples that follow the joint probability~$f_{\bX}(\bx;\mathbf{d})\mathrm{d}\mathbf{x}$ via quasi-Monte Carlo sampling or Latin hypercube sampling. Next, we construct the input-output data set from the expensive MFiX combustion simulations. In the algorithm, we assume $M<L_{S,N,m}$ data samples due to the limited computational budget. The next step is to transform the input vector $\bX$ to $\bZ$ to avoid updating the PDD's basis functions. We then use this transformation to construct the input-output data set $\{\mathbf{z}^{(l)},h(\mathbf{z}^{(l)};\mathbf{r})\}_{l=1}^M$. The algorithm proceeds to train the PDD using the input-output data via the sD-MORPH regression. We estimate the mean and variance of the output using analytical formulations \eqref{mean} and~\eqref{variance}, respectively, both derived from the PDD's coefficients. Subsequently, we conduct a derivative-free optimization, such as the Nelder-Mead method used in this work, and find a converged solution $\mathbf{d}_*^{\{i\}}$ through a single-pass training process that uses the PDD trained at the initial design $\bd_0$. Finally, the converged design vector is set as the optimal design, i.e., $\mathbf{d}^*=\mathbf{d}_*^{\{i\}}$.

The proposed RDO method achieves a significant reduction in computational cost by training the PDD surrogate model using fewer samples initially, and without requiring additional samples during the optimization. We can determine the sample size based on the computational budget, which can be lower than $L_{S,N,m}$, the number of PDD coefficients.

Other surrogate models, such as Gaussian process (GP) regression or Kriging, cannot easily replace the chosen PDD surrogate model, as we explain below. 
First, the variable transformation in step 4 of Algorithm~\ref{alg:RDO} avoids updating the PDD's basis functions when the design variable changes in every iteration. In contrast, in a GP regression or Kriging model, kernel parameters need to be updated in every iteration.
Second, the training process for GP regression model is different from that of PDD thus cannot take advantage of the computational savings by the proposed sD-MORPH method in steps 6 and 12 of Algorithm~\ref{alg:RDO}. Specifically, the proposed method finds the coefficients to PDD basis functions using limited training data by regularizing the prediction error, while GP regression aims to maximize the likelihood with respect to the kernel parameters. 
Third, the mean and standard deviation of the QoI can be obtained using the coefficients of the trained PDD. Using a GP regression model, additional samples are needed to calculate the mean and standard deviation.
In summary, the proposed method is tailored to PDD to achieve maximal computation efficiency, which may be diminished when other surrogates are used.
%
\section{Numerical results} \label{section:numerical_results}
%
In this section, we present the specifications of the design problem to maximize the thermal energy and minimize its variance for a char combustion process described in Section~\ref{section:char_combustion}, which is solved using the proposed RDO method detailed in Section~\ref{section:complete_algorithm}. For comparison, we provide results obtained via a RDO with a LASSO regression-based PDD surrogate. Section~\ref{section:numerical_setup} outlines the settings of the PDD surrogate model to approximate the thermal energy of the char combustion model. The section then describes the design variables and the random variables in detail. Finally, Section~\ref{section:results} presents the results for the RDO problem obtained using the proposed method.      
\subsection{Numerical setup}\label{section:numerical_setup}
We use the same numerical setup for PDD as that in \cite{lee2024global}, i.e., $S=2, N = 5$, and $m = 11$.  Each MFiX run to obtain a sample of the input-output data takes 425 CPU hours. We obtained 200 samples to train the PDD surrogate model, which took roughly 85,000 CPU hours. Given the limited computational budget, which determines the input-output dat set, we select the PDD's order ($m=11$) via a convergence test. This test assesses the differences to have less than 1\% in mean and variance estimates as the PDD order $m$ is incrementally increased.  We keep the ratio of the sample number to the number of unknown coefficients above 30\%, heuristically determined on Examples~1 and 2 in \cite{lee2024global}.

Given those PDD surrogate model parameters the PDD surrogate model includes $L_{S,N,m} = 606$ basis functions. The surrogate setting with $M=200$ samples results in an underdetermined regression problem~\eqref{linear} to obtain the PDD coefficients, as the sample number is only 33\% of the 606 total unknown coefficients (i.e., $M<L_{S,N,m}$). Therefore, we use sD-MORPH to compute the PDD coefficients and compare the optimization results with coefficients calculated using LASSO regression.

Table~\ref{tab:random_inputs} shows the properties of the five random inputs, namely, the height of freeboard ($X_1$), air inflow ($X_2$), the diameter of the glass bead particle ($X_3$), the diameter of the char particle ($X_4$), and char mass inflow ($X_5$), for the char combustion model. For the truncated normal random variables $X_2$ and $X_3$, the coefficients of variation are both 10\%.
\begin{table}[htbp]
   \scriptsize
   \centering
   \begin{tabular}{c|c|c|c|c|c}
    \makecell{\textbf{Random}\\\textbf{Variable}} & \textbf{Physical Quantity} & \textbf{Mean} & \makecell{\textbf{Lower}\\ \textbf{Boundary}} & \makecell{\textbf{Upper}\\ \textbf{Boundary}} & \makecell{\textbf{Probability}\\ \textbf{Distribution}} \\
    \hline
    $X_1$ & Height of freeboard (m) & 0.12  & 0.10 & 0.14 & Uniform \\
    \hline
    $X_2$ & Air inflow (m/s) & 0.825  & 0.425 & 1.225 & Truncated normal \\
    \hline
    $X_3$ & Diameter of the glass bead particle (m) & $8.0 \times 10^{-4}$  & $2.0 \times 10^{-4}$ & $1.4 \times 10^{-3}$ & Uniform \\
    \hline
    $X_4$ & Diameter of the char particle (m) & $1.0 \times 10^{-3}$ & $5.0 \times 10^{-4}$ & $1.5 \times 10^{-3}$ & Uniform \\
    \hline
    $X_5$ & Char mass inflow (kg/s) & $7.35 \times 10^{-6}$  & $1.35 \times 10^{-6}$ & $1.35 \times 10^{-5}$ & Truncated normal \\
    \end{tabular}
    \caption{Properties of the random inputs in a fluidized bed model for char combustion (reproduced from~\cite{lee2024global}).}
    \label{tab:random_inputs}
\end{table}
In the previous work~\cite{lee2024global}, we identified the air inflow ($X_2$) and the diameter of glass bead particle ($X_3$) as the most influential factors, accounting for over 90\% of the variance in the QoI (in terms of total effect Sobol' indices determined by variance-based global sensitivity analysis). In this work, we thus select the mean values of these two random variables as design variables, i.e., $d_1=\Exp[X_2]$ and $d_2=\Exp[X_3]$, where $(d_1,d_2)^{\intercal}\in \cD = [0.625, 1.025] \times [5\times 10^{-4}, 1.1\times 10^{-3}]$.

The goal of the study is to identify the two design values $(d_1, d_2)^{\intercal}$ that maximize the mean of the thermal energy $\Exp_{\bd}[Q(\bX;\bd)]$ while simultaneously minimizing its variance $\mathbb{V}\text{ar}_{\bd}[Q(\bX;\bd)]$, and given certain weights for both objectives. We solve the design problem~\eqref{eq:RDO_formulation} to obtain an optimal design $\bd^{*}=(d_1^{*},d_2^{*})^{\intercal}$, such that 
    \begin{equation}\label{opt}
    \begin{aligned}
        \bd^{*} &= \argmin_{\bd \in \cD } \ w_1\frac{\mu_0}{\mathbb{E}_{\bd}[Q(\bX;\bd)]} + w_2\frac{\sqrt{\mathbb{V}\text{ar}_{\bd}[Q(\bX;\bd)]}}{\sigma_0}
    \end{aligned}
    \end{equation}
where $\mu_0$ and $\sigma_0$ are chosen as the mean and standard deviation of the thermal energy $Q(\bX;\bd)$ at the initial design $\bd_0$.

We set the initial design vector $\bd_0=(0.825~\mathrm{m/s}, 8.0\times 10^{-4}~\mathrm{m})^{\intercal}$. To demonstrate the robustness of the proposed method with respect to sample number $M$, we tested four different sizes of training data (150, 170, 190, 200) to compute the PDD coefficients. These coefficients computed in each design iteration are used to estimate the mean and variance of the QoI, as defined in Section~\ref{section:QoI}. Since the results showed a clear quantitative trend, we herein only report the cases of 150 and 200 samples used for training to ease readability in Sections~\ref{section:case_a}~to~\ref{section:case_c}. For comparison, we also apply the process using LASSO regression in place of sD-MORPH regression. We obtain Pareto solutions for the RDO problem by considering five different combinations of weighting factors $w_1$ and $w_2$ used in the objective function. In the sD-MORPH setting, the iteration number $i$ of $\breve{\bc}^{(i)}$ in~\eqref{sdmorph} is set to at least 10 and $\lambda$ in~\eqref{new_costf} is set to 0.2.

The proposed method offers significant computational savings. High-fidelity-based RDO is computationally prohibitive, as it requires many queries of the MFiX model, each taking about 425 CPU hours. In comparison, in the proposed method, the training of the PDD surrogate using the sD-MORPH regression and the prediction using the trained PDD take 5.3 seconds and 5.5 milliseconds, respectively.

\subsection{Results}\label{section:results}

We compare the proposed method with an RDO using a LASSO regression-based PDD surrogate for five different weights $w_1$ and $w_2$ in the objective function~\eqref{opt}. The results for these five cases are presented and discussed in Sections~\ref{section:case_a}~to~\ref{section:case_de}. Section~\ref{section:Pareto} shows the Pareto optimal solutions of the LASSO-based and sD-MORPH-based RDO methods.

The accuracy of the PDD surrogate model is evaluated by using the coefficients of determination, $R^2 = 1 - SS_{\text{res}}/SS_{\text{tot}}$, where $SS_{\text{res}}$ is the residual sum of squares and $SS_{\text{tot}}$ is the total sum of squares. The $R^2$ score measures the proportion of variance in the model output, and an $R^2$ value close to one indicates a high accuracy. With 200 samples, the PDD trained using LASSO and sD-MORPH have values of $R^2 = 0.9767$ and $R^2 = 0.9967$, respectively.
\begin{figure}[ht]
    \small
        \captionsetup[subfigure]{font=small, labelformat=empty}
        \subcaptionbox{\label{fig:w_1_w_2_100_000}}[.5\textwidth]{%
\begin{tikzpicture}

\definecolor{dimgray85}{RGB}{85,85,85}
\definecolor{gainsboro229}{RGB}{229,229,229}
\definecolor{green}{RGB}{0,128,0}
\definecolor{lightgray204}{RGB}{204,204,204}
\definecolor{orange}{RGB}{255,165,0}
\definecolor{purple}{RGB}{128,0,128}

\begin{axis}[
width=0.8\linewidth,
height=0.5\linewidth,
scale only axis,
axis line style={gainsboro229},
legend cell align={left},
legend style={fill opacity=0.8, draw opacity=1, text opacity=1, draw=lightgray204, legend pos=outer north east},
tick align=outside,
tick pos=left,
x grid style={white},
xlabel=\textcolor{dimgray85}{Iteration number},
xmajorgrids,
xmin=-0.4, xmax=8.4,
xtick style={color=dimgray85},
y grid style={white},
ylabel=\textcolor{dimgray85}{Objective function value},
ymajorgrids,
ytick style={color=dimgray85}
]
\addplot [semithick, blue, mark=*, mark size=2, mark options={solid}]
table {%
0 1
1 0.91897433633086
2 0.82077166114361
3 0.673480529358237
4 0.673480529358237
5 0.673480529358237
6 0.673480529358237
7 0.673480529358237
};
\addlegendentry{sD-MORPH, 200 samples}
\addplot [semithick, purple, mark=diamond*, mark size=3, mark options={solid}]
table {%
0 1
1 0.918387609437027
2 0.822390277672692
3 0.654798003911196
4 0.654798003911196
5 0.654798003911196
6 0.654798003911196
7 0.654798003911196
};
\addlegendentry{sD-MORPH, 150 samples}
\addplot [semithick, blue, dashed, mark=*, mark size=2, mark options={solid}, line width=1.2pt]
table {%
0 1
1 0.918959987462403
2 0.825779051989323
3 0.7386851917293
4 0.7386851917293
5 0.7386851917293
6 0.7386851917293
7 0.7386851917293
};
\addlegendentry{LASSO, 200 samples}
\addplot [semithick, purple, dotted, mark=diamond*, mark size=3, mark options={solid}]
table {%
0 1
1 0.919026515105963
2 0.825713175128322
3 0.739562842272859
4 0.739562842272859
5 0.739562842272859
6 0.739562842272859
7 0.739562842272859
};
\addlegendentry{LASSO, 150 samples}
\node[anchor=north east, fill=white, draw] at (rel axis cs:0.95, 0.95) {(a) $w_1 = 1, w_2 = 0$};
\end{axis}

\end{tikzpicture}}
        \par\bigskip
        \subcaptionbox{\label{fig:w_1_w_2_000_100}}[.5\textwidth]{
\begin{tikzpicture}

\definecolor{dimgray85}{RGB}{85,85,85}
\definecolor{gainsboro229}{RGB}{229,229,229}
\definecolor{green}{RGB}{0,128,0}
\definecolor{lightgray204}{RGB}{204,204,204}
\definecolor{orange}{RGB}{255,165,0}
\definecolor{purple}{RGB}{128,0,128}

\begin{axis}[
width=0.8\linewidth,
height=0.5\linewidth,
scale only axis,
axis line style={gainsboro229},
legend cell align={left},
legend style={fill opacity=0.8, draw opacity=1, text opacity=1, draw=lightgray204},
tick align=outside,
tick pos=left,
x grid style={white},
xlabel=\textcolor{dimgray85}{Iteration number},
xmajorgrids,
xmin=-0.5, xmax=10.5,
xtick style={color=dimgray85},
y grid style={white},
ylabel=\textcolor{dimgray85}{Objective function value},
ymajorgrids,
ymin=0.743733621660056, ymax=1.01220316087333,
ytick style={color=dimgray85}
]
\addplot [semithick, blue, mark=*, mark size=2, mark options={solid}]
table {%
0 1
1 0.978552182876349
2 0.933886214639752
3 0.849898144746357
4 0.849898144746357
5 0.849898144746357
6 0.849898144746357
};
\addplot [semithick, purple, mark=diamond*, mark size=3, mark options={solid}]
table {%
0 1
1 0.96600362245729
2 0.917670739762774
3 0.831486333449609
4 0.831486333449609
5 0.831486333449609
6 0.831486333449609
};
\addplot [semithick, blue, dashed, mark=*, mark size=2, mark options={solid}]
table {%
0 1
1 0.946920731300058
2 0.946920731300058
3 0.903215815585215
4 0.812742475456799
5 0.774627972086897
6 0.766888931567799
7 0.766888931567799
8 0.763340324609246
9 0.763340324609246
10 0.763340324609246
};
\addplot [semithick, purple, dotted, mark=diamond*, mark size=3, mark options={solid}]
table {%
0 1
1 0.945095467839414
2 0.945095467839414
3 0.899599346502056
4 0.807240433245689
5 0.773326576078385
6 0.756868245221082
7 0.756868245221082
8 0.755936782533387
9 0.755936782533387
10 0.755936782533387
};
\node[anchor=north east, fill=white, draw] at (rel axis cs:0.95, 0.95) {(b) $w_1 = 0, w_2 = 1$};
\end{axis}

\end{tikzpicture}}
        \subcaptionbox{\label{fig:w_1_w_2_050_050}}[.5\textwidth]{
\begin{tikzpicture}

\definecolor{dimgray85}{RGB}{85,85,85}
\definecolor{gainsboro229}{RGB}{229,229,229}
\definecolor{green}{RGB}{0,128,0}
\definecolor{lightgray204}{RGB}{204,204,204}
\definecolor{orange}{RGB}{255,165,0}
\definecolor{purple}{RGB}{128,0,128}

\begin{axis}[
width=0.8\linewidth,
height=0.5\linewidth,
scale only axis,
axis line style={gainsboro229},
tick align=outside,
tick pos=left,
x grid style={white},
xlabel=\textcolor{dimgray85}{Iteration number},
xmajorgrids,
xmin=-0.4, xmax=8.4,
xtick style={color=dimgray85},
y grid style={white},
ylabel=\textcolor{dimgray85}{Objective function value},
ymajorgrids,
ymin=0.970717621706061, ymax=1.00139441351934,
ytick style={color=dimgray85}
]
\addplot [semithick, blue, mark=*, mark size=2, mark options={solid}]
table {%
0 1
1 0.988686602295958
2 0.988686602295958
3 0.984368536436484
4 0.980570273868798
5 0.980570273868798
6 0.977096565447862
7 0.977096565447862
8 0.977096565447862
};
\addplot [semithick, purple, mark=diamond*, mark size=3, mark options={solid}]
table {%
0 1
1 0.983058712810876
2 0.981349186927202
3 0.978415242085067
4 0.978415242085067
5 0.972112021333937
6 0.972112021333937
};
\addplot [semithick, blue, dashed, mark=*, mark size=2, mark options={solid}]
table {%
0 1
1 0.996264238966663
2 0.996264238966663
3 0.996264238966663
4 0.995555316990469
5 0.995555316990469
6 0.995555316990469
7 0.995555316990469
};
\addplot [semithick, purple, dotted, mark=diamond*, mark size=3, mark options={solid}]
table {%
0 1
1 0.99642859564928
2 0.994263544054497
3 0.994263544054497
4 0.993331321178681
5 0.993100304772822
6 0.993100304772822
7 0.993100304772822
};
\node[anchor=north east, fill=white, draw] at (rel axis cs:0.95, 0.65) {(c) $w_1 = 0.5, w_2 = 0.5$};
\end{axis}

\end{tikzpicture}}
        \caption{Change in objective function value with increasing iterations for different combinations of $w_1$ and $w_2$ in \eqref{opt}, comparison of LASSO-based and sD-MORPH-based design methods over different sample numbers (150 and 200).} 
\end{figure}
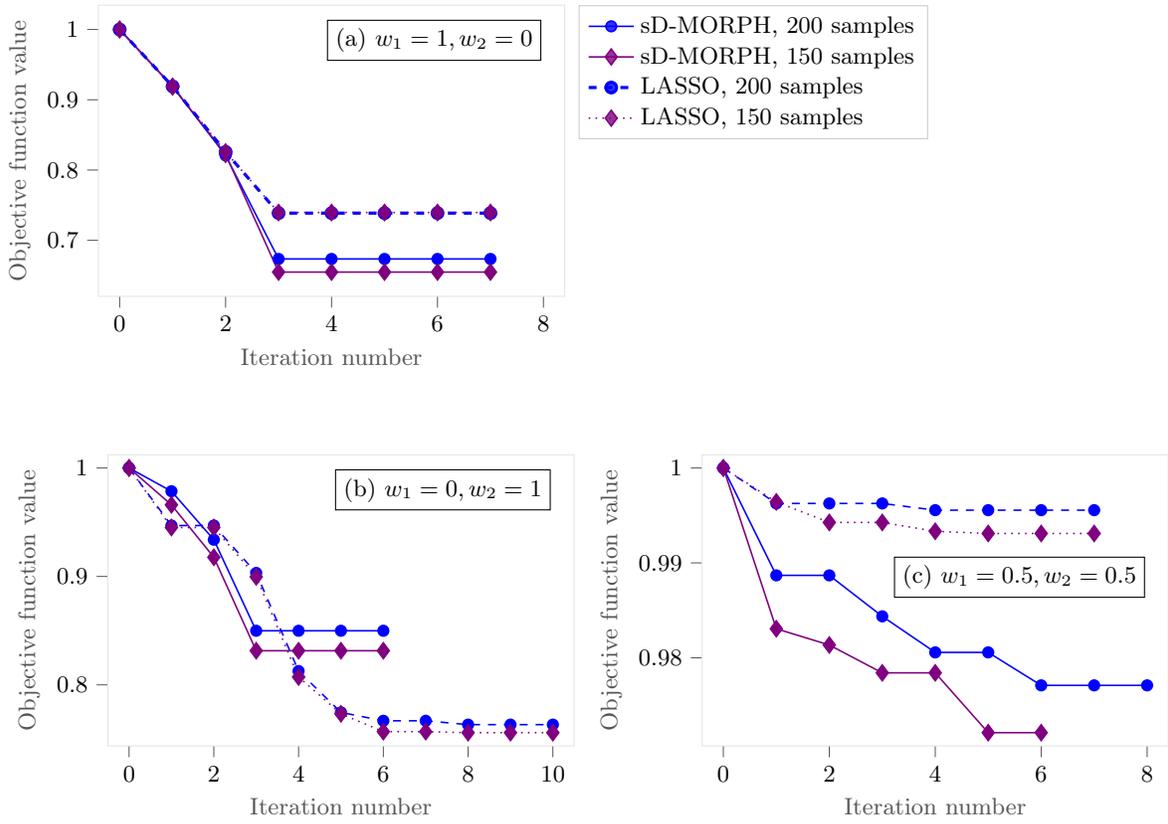
\subsubsection{Case (a): $w_1=1$ and $w_2 = 0$} \label{section:case_a}
This combination of weighting factors maximizes the mean of thermal energy, without considering its standard deviation. Figure~\ref{fig:w_1_w_2_100_000} compares the optimization results obtained using LASSO-based and sD-MORPH-based RDO methods. We observe that the objective function value consistently decreases across all cases with different regression methods and different numbers of samples. 

Using 200 training samples, at the initial design vector $\bd_0 = (0.825~\mathrm{m/s},~ 8.0\times 10^{-4}~\mathrm{m})^{\intercal}$, the mean values of the thermal energy esimated with the PDD surrogate coefficients using LASSO and sD-MORPH regressions are 1437.56J and 1436.78J, respectively. The optimal design vector is $\bd^* = (1.025~\mathrm{m/s},~5.0\times 10^{-4}~\mathrm{m})^{\intercal}$ for both methods, where the corresponding mean values of the thermal energy are 2017.08J (40.31\% increase from initial design) and 2133.37J (48.48\% increase), respectively. These results, shown in Table~\ref{table:results_w_1_1}, indicate that the sD-MORPH-based RDO method increases the mean of the thermal energy more. The results of RDO methods using 150 training samples are similar.
    \begin{table}[ht]
        \centering
        \begin{tabular}{c|c|c}
        \hline
                 & LASSO & sD-MORPH \\ \hline
        Initial design $\bd_0 = (0.825~\mathrm{m/s},~ 8.0\times 10^{-4}~\mathrm{m})^{\intercal}$ & 1437.56 J   & 1436.78 J   \\ \hline
        Optimal design $\bd^* = (1.025~\mathrm{m/s},~5.0\times 10^{-4}~\mathrm{m})^{\intercal}$  & 2017.08 J  & 2133.37 J   \\ \hline
        Percentage change at $\bd^*$ compared to at $\bd_0$ & +40.31\% & + 48.48\% \\ \hline
        \end{tabular}
        \caption{Case (a): mean values of thermal energy at the initial design and the optimal design obtained using PDD coefficients trained with 200 samples using LASSO and sD-MORPH regressions.}
        \label{table:results_w_1_1}
    \end{table}

Note that the mean values of thermal energy are estimated with the PDD coefficients calculated using LASSO or sD-MORPH regressions, see~\eqref{mean}. Therefore, even with the same design vector, the values may be different.
\subsubsection{Case (b): $w_1=0$ and $w_2 = 1$}\label{section:case_b}

This combination of weighting factors minimizes the standard deviation of thermal energy, without considering its mean. With 200 training samples, the standard deviations of the thermal energy obtained using the LASSO-based and sD-MORPH-based RDO methods at the initial design are 185.13J and 198.67J, respectively. The optimal design vectors for the two methods are $(0.625~\mathrm{m/s},~8.25\times 10^{-4}~\mathrm{m})^{\intercal}$ and $(0.709~\mathrm{m/s},~7.96\times 10^{-4}~\mathrm{m})^{\intercal}$, respectively; the corresponding standard deviations of the thermal energy are 141.31J (23.67\% decrease) and 168.85J (15.01\% decrease), respectively. These results are shown in Table~\ref{table:results_w_2_1}. Compared to the LASSO-based method, the proposed sD-MORPH based method results in about 8\% less reduction in the objective function value. This is because LASSO is limited by the number of non-zero values which is constrained by the number of training samples. As a result, the standard deviation (estimated by \eqref{variance}) is often underestimated. In contrast, the sD-MORPH solution can produce small non-zero values, providing a more accurate representation of the true values. 
\begin{table}[ht]
    \centering
    \begin{tabular}{c|c}
    \hline
    \multicolumn{2}{c}{LASSO} \\ \hline
    $\bd_0 = (0.825~\mathrm{m/s},~ 8.0\times 10^{-4}~\mathrm{m})^{\intercal}$ & $\sqrt{\mathbb{V}\text{ar}[Q]}$ = 185.13J \\ \hline
    $\bd^* = (0.625~\mathrm{m/s},~8.25\times 10^{-4}~\mathrm{m})^{\intercal}$ & $\sqrt{\mathbb{V}\text{ar}[Q]}$ = 141.31J \\ \hline
    Percentage change at $\bd^*$ compared to at $\bd_0$ & -23.67\% \\ \hline \hline
    \multicolumn{2}{c}{sD-MORPH} \\ \hline
    
    $\bd_0 = (0.825~\mathrm{m/s},~ 8.0\times 10^{-4}~\mathrm{m})^{\intercal}$ & $\sqrt{\mathbb{V}\text{ar}[Q]}$ = 198.67J \\ \hline
    $\bd^* = (0.709~\mathrm{m/s},~7.96\times 10^{-4}~\mathrm{m})^{\intercal}$ & $\sqrt{\mathbb{V}\text{ar}[Q]}$ = 168.85J \\ \hline
    Percentage change at $\bd^*$ compared to at $\bd_0$ & -15.01\% \\ \hline
    \end{tabular}
    \caption{Case (b): standard deviations of the thermal energy at initial design and optimal designs obtained using PDD coefficients trained with 200 samples using LASSO and sD-MORPH regressions.}
    \label{table:results_w_2_1}
\end{table}

In Figure~\ref{fig:w_1_w_2_000_100},  we compare the optimization results obtained using LASSO and sD-MORPH-based RDO methods. The sD-MORPH version shows a faster drop of about 15\% in the objective function within three iterations, whereas the LASSO version reduces the cost function by about 10\% within the first three iterations. By the sixth iteration, the sD-MORPH version has converged, while the LASSO version is approaching convergence by the tenth iteration. The results of RDO methods with 150 training samples are similar.

\subsubsection{Case (c): $w_1=0.5$ and $w_2 = 0.5$} \label{section:case_c}

This combination of weighting parameters corresponds to an optimization where the mean and standard deviation of the thermal energy are equally prioritized. Figure~\ref{fig:w_1_w_2_050_050} shows the comparison of optimizations with LASSO-based and sD-MORPH-based methods.
The objective function value consistently decreases across all cases with different regression methods and different numbers of samples. 

For the LASSO-based method with 200 samples, the optimal design vector is $\bd^* = (0.816~\mathrm{m/s}, 8.24\times 10^{-4}~\mathrm{m})^{\intercal}$. The corresponding standard deviation is 184.32J, almost the same as 185.13J at the initial design. By comparison, for sD-MORPH-based method with 200 samples, the optimal design vector is $\bd^* = (0.768~\mathrm{m/s}, 8.10\times 10^{-4}~\mathrm{m})^{\intercal}$. The corresponding standard deviation is 173.94J (12.45\% decrease). This shows that the optimization for equally weighted mean and standard deviation with sD-MORPH regression effectively reduces the standard deviation. The mean values of the thermal energy corresponding to the LASSO-based and sD-MORPH-based RDO methods are 1437.56J and 1436.55J at the initial design, respectively; at the optimal designs, they are 1417.83J (1.37\% decrease) and 1336.28J (6.98\% decrease), respectively. These results are shown in Table~\ref{table:results_w_1_0.5}. The results when using 150 samples are similar.

RDO minimizes the objective function, in our case the sum of the normalized mean and variance. With 200 samples, the objective function value decreases from 1 at the initial design to 0.9956 and 0.9771 at the optimal designs for LASSO and sD-MORPH cases, respectively. The goal of RDO is to minimize the sensitivity of design performance to variations and to find designs that perform consistently across different conditions; this does not guarantee a decrease in standard deviation and an increase in mean. In this case, both mean and standard deviation decrease after optimization, however, the objective function value decreases, which meets the optimization objective.
\begin{table}[ht]
    \centering
    \begin{tabular}{c|c|c}
    \hline
    \multicolumn{3}{c}{LASSO} \\ \hline
    $\bd_0 = (0.825~\mathrm{m/s},~ 8.0\times 10^{-4}~\mathrm{m})^{\intercal}$ & $\sqrt{\mathbb{V}\text{ar}[Q]}$ = 185.13J & $\mathbb{E}[Q]$ = 1437.56J \\ \hline
    $\bd^* = (0.816~\mathrm{m/s}, 8.24\times 10^{-4}~\mathrm{m})^{\intercal}$ & $\sqrt{\mathbb{V}\text{ar}[Q]}$ = 184.32J & $\mathbb{E}[Q]$ = 1417.83J\\ \hline
    Percentage change at $\bd^*$ compared to at $\bd_0$ & -0.44\% & -1.37\%\\ \hline \hline
    \multicolumn{3}{c}{sD-MORPH} \\ \hline
    
    $\bd_0 = (0.825~\mathrm{m/s},~ 8.0\times 10^{-4}~\mathrm{m})^{\intercal}$ & $\sqrt{\mathbb{V}\text{ar}[Q]}$ = 198.67J & $\mathbb{E}[Q]$ = 1436.55J\\ \hline
    $\bd^* = (0.768~\mathrm{m/s}, 8.10\times 10^{-4}~\mathrm{m})^{\intercal}$ & $\sqrt{\mathbb{V}\text{ar}[Q]}$ = 173.94J & $\mathbb{E}[Q]$ = 1336.28J\\ \hline
    Percentage change at $\bd^*$ compared to at $\bd_0$ & -12.45\% & -6.98\%\\ \hline
    \end{tabular}
    \caption{Case (c): means and standard deviations of the thermal energy at initial design and optimal designs obtained using PDD coefficients trained with 200 samples using LASSO and sD-MORPH regressions.}
    \label{table:results_w_1_0.5}
\end{table}

\subsubsection{Other cases} \label{section:case_de}
Results for other two cases ($w_1=0.25$, $w_2 = 0.75$; $w_1=0.75$, $w_2 = 0.25$) are shown in Appendix~\ref{section:appendix_1}.

\subsubsection{Numerical verification}
To verify the surrogate-based optimization results with the high-fidelity PIC solver, we calculate unbiased reference values for mean and standard deviation. For the initial design $\bd_0 = (0.825~\mathrm{m/s},~ 8.0\times 10^{-4}~\mathrm{m})^{\intercal}$, we evaluate 200 Monte Carlo simulation results, where the initial mean is 1414.64J and the variance 195.24J.  Comparing with the values in Table~\ref{table:results_w_1_0.5} we observe that the mean and standard deviation estimates from sD-MORPH are closer to the reference values than estimations by the LASSO-based surrogate.

For case (c), we perform a numerical verification by conducting additional high-fidelity PIC simulations using the optimal designs of both the LASSO-based and sD-MORPH-based RDO. For the optimal design $\bd^* = (0.816~\mathrm{m/s}, 8.24\times 10^{-4}~\mathrm{m})^{\intercal}$ obtained by the LASSO-based method, a Monte Carlo simulation with 199 samples produced a mean of 1404.93J and a standard deviation of 191.64J, yielding an objective function value of 1.0062. In contrast, for the optimal design $\bd^* = (0.768~\mathrm{m/s}, 8.10\times 10^{-4}~\mathrm{m})^{\intercal}$ obtained by the sD-MORPH-based method, the Monte Carlo simulation provided a mean of 1322.51J and a standard deviation of 179.83J; the corresponding objective function value is 0.9129.
Compared to the values in Table~\ref{table:results_w_1_0.5}, the LASSO-based method estimates the mean as 0.92\% higher and the standard deviation as 3.82\% lower compared to the Monte Carlo simulation. In contrast, the sD-MOPRH-based method estimates are 1.03\% higher for the mean and 3.38\% lower for the standard deviation relative to Monte Carlo simulation. Despite a low number of Monte Carlo samples, this shows that: (1) sD-MORPH-based method is more accurate in terms of estimating the standard deviation value; and (2) sD-MORPH-based method effectively reduces the standard deviation and the objective function value.

\subsubsection{Pareto solutions} \label{section:Pareto}
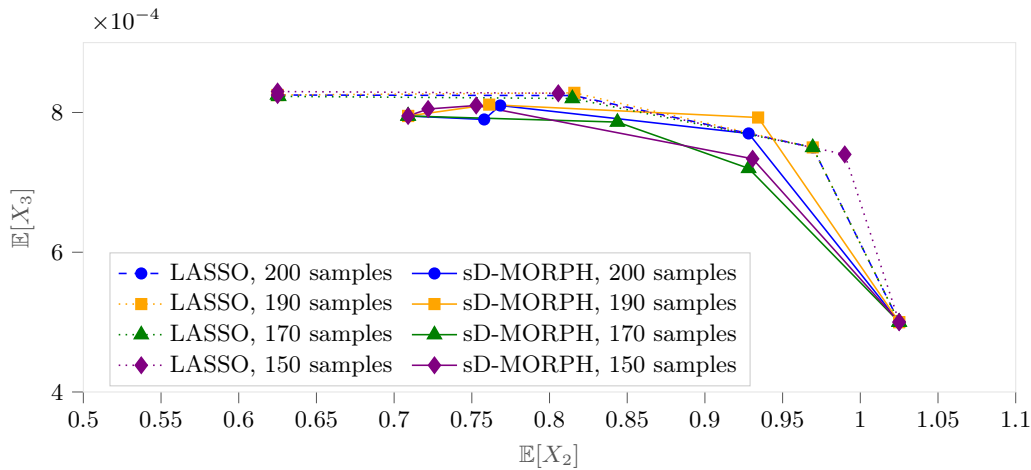
\begin{figure}
        \centering
        \small
\begin{tikzpicture}

\definecolor{dimgray85}{RGB}{85,85,85}
\definecolor{gainsboro229}{RGB}{229,229,229}
\definecolor{green01270}{RGB}{0,127,0}
\definecolor{lightgray204}{RGB}{204,204,204}
\definecolor{orange}{RGB}{255,165,0}
\definecolor{purple}{RGB}{128,0,128}

\begin{axis}[
width=0.8\linewidth,
height=0.3\linewidth,
at={(0\fheight,0\fheight)},
scale only axis,
axis line style={gainsboro229},
legend columns=2,
legend cell align={left},
legend style={draw opacity=1, text opacity=1, draw=lightgray204, /tikz/column 2/.style={column sep=5pt}, {at={(rel axis cs:0.028, 0.016)}, anchor=south west}},
tick align=outside,
tick pos=left,
tick scale binop=\times,
x grid style={white},
xlabel=\textcolor{dimgray85}{\(\displaystyle \Exp[X_2]\)},
xmajorgrids,
xmin=0.5, xmax=1.1,
xtick style={color=dimgray85},
y grid style={white},
ylabel=\textcolor{dimgray85}{\(\displaystyle \Exp[X_3]\)},
ymajorgrids,
ymin=0.0004, ymax=0.0009,
ytick style={color=dimgray85}
]
\addplot [semithick, blue, dashed, mark=*, mark size=2, mark options={solid}]
table {%
0.625 0.000824999999999999
0.625 0.000824999999999999
0.815654296875 0.00082421875
0.969375 0.00075
1.025 0.0005
};
\addlegendentry{LASSO, 200 samples}
\addplot [semithick, blue, mark=*, mark size=2, mark options={solid}]
table {%
0.708984375 0.000795000000000001
0.75796875 0.000790000000000001
0.76828125 0.00081
0.928125 0.00077
1.025 0.0005
};
\addlegendentry{sD-MORPH, 200 samples}
\addplot [semithick, orange, dotted, mark=square*, mark size=2, mark options={solid}]
table {%
0.625 0.000824999999999999
0.625 0.000824999999999999
0.8159765625 0.000828125
0.969375 0.00075
1.025 0.0005
};
\addlegendentry{LASSO, 190 samples}
\addplot [semithick, orange, mark=square*, mark size=2, mark options={solid}]
table {%
0.708984375 0.000795000000000001
0.708984375 0.000795000000000001
0.76119140625 0.00081125
0.934248046875 0.0007928125
1.025 0.0005
};
\addlegendentry{sD-MORPH, 190 samples}
\addplot [semithick, green01270, dotted, mark=triangle*, mark size=3, mark options={solid}]
table {%
0.625 0.000824999999999999
0.625674438476562 0.0008230615234375
0.8146875 0.00082
0.969375 0.00075
1.025 0.0005
};
\addlegendentry{LASSO, 170 samples}
\addplot [semithick, green01270, mark=triangle*, mark size=3, mark options={solid}]
table {%
0.708984375 0.000795000000000001
0.708984375 0.000795000000000001
0.84369140625 0.00078625
0.928125 0.00072
1.025 0.0005
};
\addlegendentry{sD-MORPH, 170 samples}
\addplot [semithick, purple, dotted, mark=diamond*, mark size=3, mark options={solid}]
table {%
0.625 0.000824999999999999
0.625 0.00083
0.8056640625 0.0008275
0.99 0.00074
1.025 0.0005
};
\addlegendentry{LASSO, 150 samples}
\addplot [semithick, purple, mark=diamond*, mark size=3, mark options={solid}]
table {%
0.708984375 0.000795000000000001
0.721875 0.000805
0.7528125 0.00081
0.930703125 0.00073375
1.025 0.0005
};
\addlegendentry{sD-MORPH, 150 samples}
\end{axis}

\end{tikzpicture}
        \caption{Pareto solutions in the design space $\mathbb{E}[X_2]=d_1$ and $\mathbb{E}[X_3]=d_2$ in comparison of sD-MORPH-and LASSO-based design methods over different sample numbers (150, 170, 190, 200).}
        \label{fig:Pareto}
    \end{figure}

Figure~\ref{fig:Pareto} shows Pareto optimal solutions obtained using sD-MORPH-based and LASSO-based RDO methods, each with different sample numbers. Pareto solutions obtained with sD-MORPH-based RDO method fall within that of LASSO-based RDO method, indicating that sD-MORPH exhibits less variation across different sample sizes. This demonstrates the robustness of sD-MORPH in maintaining optimal performance regardless of the number of samples. The Pareto solutions provide guidance on how to adjust the operation parameters ($X_2$, air inflow and $X_3$, diameter of the glass bead particle) in order to achieve consistently high thermal energy output while minimizing its variation.

\section{Conclusion and future work} \label{section:conclusion}
We proposed a robust design optimization (RDO) problem for a char combustion process and solved it in a limited-data setting with a suitable surrogate model. We used an sD-MORPH regression to calculate the coefficients of the PDD surrogate model in the optimization and reformulated the optimization problem using transformed input random variables combined with a single-pass training process. 
We obtained Pareto solutions with different weights for the mean and standard deviation, which helps to choose operation parameters to produce a consistently high thermal energy from char combustion. We also conducted comparisons between RDO methods that use either a LASSO or an sD-MORPH-based surrogate model. The LASSO-based and sD-MORPH-based RDO yielded comparable improvements in thermal energy production when maximizing only the mean without accounting for the standard deviation. In another test where we only minimized the standard deviation, the sD-MORPH-based RDO converged faster than the LASSO-based approach. In a third case, we optimized both the mean and standard deviation of thermal energy, and found that the sD-MORPH-based method more effectively reduces standard deviation. We also found that the sD-MORPH surrogate more accurately estimates both the mean and standard deviation, as confirmed by some verification simulations.


In the proposed method, we train the PDD surrogate model once to represent the output over the entire design domain. This limits us to a narrow design space where the PDD can accurately represent the output while using fewer data points for RDO, as demonstrated by the optimization of the char combustion problem.
The proposed method can be scaled to solve problems with higher-dimensional input spaces or additional objectives as follows.  
For problems with a higher-dimensional input space, a larger number of training data from physics-based model simulations may be necessary as the number of sensitive input variables could be high, resulting in a larger number of basis functions. 
For optimization problems with multiple objectives and/or constraints, PDD surrogate models corresponding to the quantities in the objectives and constraints can be constructed and trained.
A promising direction is to extend this work to include constraints in the optimization, such as reliability or risk measures.
%
%
%
%
\section*{Statements and declarations}
\subsection*{Funding}\label{ack} 
Y. Guo and B. Kramer were in part financially supported by the Korea Institute for Advancement of Technology (KIAT) through the International Cooperative R\&D program (No. P0019804, Digital twin based intelligent unmanned facility inspection solutions). D. Lee was supported by the research fund of Hanyang University (HY-202300000003174).
\subsection*{CRediT authorship contribution statement}
Yulin Guo: Formal analysis, Investigation, Data Curation, Software, Visualization, Writing - original draft; Dongjin Lee: Conceptualization, Methodology, Data Curation, Software, Visualization, Writing - original draft, Funding acquisition; Boris Kramer: Conceptualization, Writing - review and editing, Funding acquisition, Project administration, Supervision. 
\subsection*{Conflict of interest}
Boris Kramer reports a relationship with ASML Holding US that includes: consulting or advisory. Other authors declare that they have no known competing financial interests or personal relationships that could have appeared to influence the work reported in this paper.
%
\subsection*{Replication of results}
A Python code file corresponding to one optimization case of the proposed method and MFiX simulation code for the char combustion example are available to download in the following Github repository: {https://github.com/yulin-g/sD-MORPH-optimization}.
\subsection*{Data availability}
Data can be replicated by using the code shared on GitHub.
\subsection*{Ethics approval and consent to participate}
Not applicable.
\begin{appendices}
\section{Heat capacity calculations} \label{section:appendix_cp}
The specific heat capacity $C_{pi}$ is calculate as $C_{pi} = \frac{R}{Mw_i}\left(A+BT+CT^2+DT^3+ET^4\right)$, where $R$ is the gas constant, $Mw_i$ is the molecular weight of the gas, $A,B,C,D$, and $E$ are coefficients obtained from \cite{burcat2005third}, and $T$ is the temperature of the gas. 

The molecular weights are 31.9988 g/mol (O\textsubscript{2}), 28.0134 g/mol (N\textsubscript{2}), 28.0104 g/mol (CO), 44.0098 g/mol (CO\textsubscript{2}), and 18.0153 g/mol (H\textsubscript{2}O). The coefficients are listed in Tables~\ref{table:cp_coeff_T_l_1000}~and~\ref{table:cp_coeff_T_g_1000}.

\begin{table}[h!]
\centering
\begin{tabular}{c|c|c|c|c|c}
\hline
Gas & $A$ & $B$ & $C$ & $D$ & $E$ \\ \hline
O\textsubscript{2}    & $3.7825$ & $-2.9970 \times 10^{-3}$ & $9.8473 \times 10^{-6}$ & $-9.6813 \times 10^{-9}$ & $3.2437 \times 10^{-12}$ \\ \hline
N\textsubscript{2}    & $3.5310$ & $-1.2366 \times 10^{-4}$ & $-5.0300 \times 10^{-7}$ & $2.4353 \times 10^{-9}$ & $-1.4088 \times 10^{-12}$ \\ \hline
CO    & $3.5795$ & $-6.1035 \times 10^{-4}$ & $1.0168 \times 10^{-6}$ & $9.0701 \times 10^{-10}$ & $-9.0442 \times 10^{-13}$ \\ \hline
CO\textsubscript{2}   & $2.3568$ & $8.9841 \times 10^{-3}$ & $-7.1221 \times 10^{-6}$ & $2.4573 \times 10^{-9}$ & $-1.4289 \times 10^{-13}$ \\ \hline
H\textsubscript{2}O   & $4.1986$ & $-2.0364 \times 10^{-3}$ & $6.5203 \times 10^{-6}$ & $-5.4879 \times 10^{-9}$ & $1.7720 \times 10^{-12}$ \\ \hline
\end{tabular}
\caption{Coefficients for calculating $C_{pi}$ when gas temperature $T<1,000.K$}\label{table:cp_coeff_T_l_1000}
\end{table}

\begin{table}[h!]
\centering
\begin{tabular}{c|c|c|c|c|c}
\hline
Gas & $A$ & $B$ & $C$ & $D$ & $E$ \\ \hline
O\textsubscript{2}    & $3.6610$ & $6.5637 \times 10^{-4}$ & $-1.4115 \times 10^{-7}$ & $2.0580 \times 10^{-11}$ & $-1.2991 \times 10^{-15}$ \\ \hline
N\textsubscript{2}    & $2.9526$ & $1.3969 \times 10^{-3}$ & $-4.9263 \times 10^{-7}$ & $7.8601 \times 10^{-11}$ & $-4.6076 \times 10^{-15}$ \\ \hline
CO    & $3.0485$ & $1.3517 \times 10^{-3}$ & $-4.8579 \times 10^{-7}$ & $7.8854 \times 10^{-11}$ & $-4.6981 \times 10^{-15}$ \\ \hline
CO\textsubscript{2}   & $4.6365$ & $2.7415 \times 10^{-3}$ & $-9.9590 \times 10^{-7}$ & $1.6039 \times 10^{-10}$ & $-9.1620 \times 10^{-15}$ \\ \hline
H\textsubscript{2}O   & $2.6770$ & $2.9732 \times 10^{-3}$ & $-7.7377 \times 10^{-7}$ & $9.4434 \times 10^{-11}$ & $-4.2690 \times 10^{-15}$ \\ \hline
\end{tabular}
\caption{Coefficients for calculating $C_{pi}$ when gas temperature $T \geq 1,000. K$}\label{table:cp_coeff_T_g_1000}
\end{table}

\section{Additional numerical results} \label{section:appendix_1}
\begin{figure}[h!]
    \small
        \captionsetup[subfigure]{font=small}
        \subcaptionbox{$w_1=0.25$, $w_2 = 0.75$ \label{fig:w_1_w_2_025_075}}[.5\textwidth]{%
\begin{tikzpicture}

\definecolor{dimgray85}{RGB}{85,85,85}
\definecolor{gainsboro229}{RGB}{229,229,229}
\definecolor{green}{RGB}{0,128,0}
\definecolor{lightgray204}{RGB}{204,204,204}
\definecolor{orange}{RGB}{255,165,0}
\definecolor{purple}{RGB}{128,0,128}

\begin{axis}[
width=0.8\linewidth,
height=0.5\linewidth,
scale only axis,
axis line style={gainsboro229},
legend cell align={left},
legend style={fill opacity=0.8, draw opacity=1, text opacity=1, draw=lightgray204, font = \scriptsize},
tick align=outside,
tick pos=left,
x grid style={white},
xlabel=\textcolor{dimgray85}{Iteration number},
xmajorgrids,
xmin=-0.5, xmax=10.5,
xtick style={color=dimgray85},
y grid style={white},
ylabel=\textcolor{dimgray85}{Objective function value},
ymajorgrids,
ymin=0.892452722358283, ymax=1.00512129893532,
ytick style={color=dimgray85}
]
\addplot [semithick, blue, mark=*, mark size=2, mark options={solid}]
table {%
0 1
1 0.983668200630836
2 0.960998635123752
3 0.928624963824011
4 0.928624963824011
5 0.928624963824011
};
\addlegendentry{sD-MORPH, 200 samples}
\addplot [semithick, purple, mark=diamond*, mark size=3, mark options={solid}]
table {%
0 1
1 0.975225551265949
2 0.975225551265949
3 0.90990455537252
4 0.90990455537252
5 0.908029030320155
6 0.908029030320155
};
\addlegendentry{sD-MORPH, 150 samples}
\addplot [semithick, blue, dashed, mark=*, mark size=2, mark options={solid}]
table {%
0 1
1 0.973660013189619
2 0.973660013189619
3 0.953677768157468
4 0.918490676306068
5 0.908302285417783
6 0.905192096740685
7 0.905192096740685
8 0.902852313306603
9 0.902852313306603
10 0.902852313306603
};
\addlegendentry{LASSO, 200 samples}
\addplot [semithick, purple, dotted, mark=diamond*, mark size=3, mark options={solid}]
table {%
0 1
1 0.972131286544106
2 0.972131286544106
3 0.950987000279997
4 0.914244283856213
5 0.907306176145865
6 0.897574021293603
7 0.897574021293603
8 0.897574021293603
9 0.897574021293603
10 0.897574021293603
};
\addlegendentry{LASSO, 150 samples}
\end{axis}

\end{tikzpicture}}%
        \subcaptionbox{$w_1=0.75$, $w_2 = 0.25$\label{fig:w_1_w_2_075_025}}[.5\textwidth]{
\begin{tikzpicture}

\definecolor{dimgray85}{RGB}{85,85,85}
\definecolor{gainsboro229}{RGB}{229,229,229}
\definecolor{green}{RGB}{0,128,0}
\definecolor{lightgray204}{RGB}{204,204,204}
\definecolor{orange}{RGB}{255,165,0}
\definecolor{purple}{RGB}{128,0,128}

\begin{axis}[
width=0.8\linewidth,
height=0.5\linewidth,
scale only axis,
axis line style={gainsboro229},
legend cell align={left},
legend columns=2,
legend style={fill opacity=0.8, draw opacity=1, text opacity=1, draw=lightgray204, /tikz/column 2/.style={column sep=5pt}},
tick align=outside,
tick pos=left,
x grid style={white},
xlabel=\textcolor{dimgray85}{Iteration number},
xmajorgrids,
xmin=-0.5, xmax=10.5,
xtick style={color=dimgray85},
y grid style={white},
ylabel=\textcolor{dimgray85}{Objective function value},
ymajorgrids,
ymin=0.912904120568169, ymax=1.00414742283009,
ytick style={color=dimgray85}
]
\addplot [semithick, blue, mark=*, mark size=2, mark options={solid}]
table {%
0 1
1 0.976548068927963
2 0.958970628162039
3 0.958970628162039
4 0.958970628162039
5 0.947161776386773
6 0.947161776386773
7 0.947161776386773
8 0.947161776386773
};
\addplot [semithick, blue, dashed, mark=*, mark size=2, mark options={solid}]
table {%
0 1
1 0.976226474208961
2 0.936102430704197
3 0.934895615224772
4 0.934895615224772
5 0.923703268547204
6 0.923703268547204
7 0.923703268547204
};
\addplot [semithick, purple, mark=diamond*, mark size=3, mark options={solid}]
table {%
0 1
1 0.965798607730279
2 0.946247570991934
3 0.946247570991934
4 0.943157613073607
5 0.943157613073607
6 0.943157613073607
7 0.943157613073607
8 0.938757744585098
9 0.938757744585098
10 0.938757744585098
};
\addplot [semithick, purple, dotted, mark=diamond*, mark size=3, mark options={solid}]
table {%
0 1
1 0.975545253204161
2 0.934790884807164
3 0.934790884807164
4 0.934790884807164
5 0.918992301653874
6 0.917051543398256
7 0.917051543398256
8 0.917051543398256
};
\end{axis}

\end{tikzpicture}}
        \caption{Change in objective function value with increasing iterations for different combinations of $w_1$ and $w_2$ in \eqref{opt}, comparison of LASSO-based and sD-MORPH-based design methods over different sample numbers (170 and 190).} 
\end{figure}
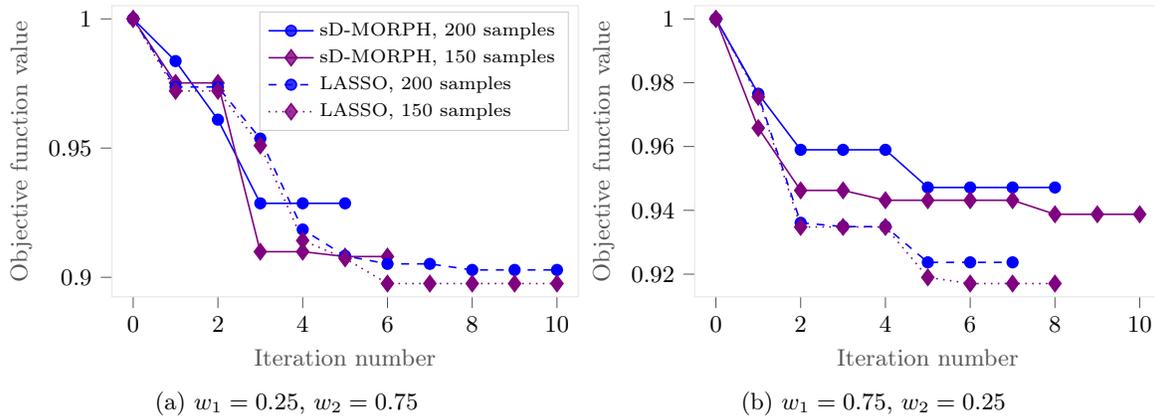

For additional weighting parameter settings ($w_1=0.25$, $w_2 = 0.75$; $w_1=0.75$, $w_2 = 0.25$), the changes of objective function value versus the iteration number are shown in Figures~\ref{fig:w_1_w_2_025_075}~and~\ref{fig:w_1_w_2_075_025}.

\section*{Nomenclature}

\renewcommand{\arraystretch}{1.3} 

\begin{longtable}{p{2cm} p{5cm} p{2.5cm} p{4.5cm}}
	$\mathbf{a}(t)$ & Potential sD-MORPH solution & $\bA^+$ & Intermediate matrix for calculating $\breve{\bc}$, details in \cite{lee2024global} \\
	$\mathbb{A}^{N}$ & Subdomain of $\real^N$ & $\mathcal{B}^{N}$ & Borel $\sigma$-field on $\mathbb{A}^{N}$ \\
	$\breve{\mathbf{c}}$ & sD-MOPRPH estimates of the PDD coefficients & $\mathbf{c}_0$ & LASSO estimates of the PDD coefficients \\
	$c$ & Expansion coefficient & $C_{\text{CO}}, C_{\text{O}_2}, C_{\text{H}_2\text{O}}$ & Mass concentrations of CO, O$_2$, H$_2$O, $\mathrm{kmole/m^3}$ \\
	$C_p$ & Specific heat capacity of the gas mixture $\mathrm{J/(kg\cdot K)}$ & $C_{pi}$ & Specific heat capacity of a gas component, $\mathrm{J/(kg\cdot K)}$ \\
	$\mathcal{D}$ & Design space & $\mathbf{d}$ & Vector of design variables \\
	$d_p$ & Particle size, $\mathrm{m}$ & $D_o$ & Oxygen-nitrogen mixture diffusion coefficient \\
	$\bar{\mathbf{E}}$ & Intermediate matrix for calculating $\breve{\bc}$, details in~\cite{lee2024global} & $\mathcal{F}_{\bd}$ & $\sigma$-field on $\Omega_{\bd}$ \\
	$\bar{\mathbf{F}}$ & Intermediate matrix for calculating $\breve{\bc}$, details in~\cite{lee2024global} & $\mathbf{g}$ & Vector of mean values \\
	$h(\cdot)$ & Thermal energy as a function of transformed random variables & $h_{S,m}$ & $S$-variate, $m$-th order PDD approximation of $h$ \\
	$\breve{\mathcal{K}}$ & Cost function in sD-MORPH regression & $L$ & Number of PDD basis functions \\
	$m$ & Number of highest PDD order & $m_{ci}$ & Unreacted char mass, $\mathrm{g}$ \\
	$\dot{m}$ & Mass flow rate, $\mathrm{g/s}$ & $M_i$ & Mole fraction of a gas component, $\mathrm{g/mol}$ \\
	$N$ & Number of random inputs & $\Omega_{\bd}$ & Sample space \\
	$p_o$ & Oxygen partial pressure, $\mathrm{Pa}$ & $\mathbb{P}_{\bd}$ & A family of probability measures \\
	$P_s$ & Pressure constant & $Q$ & Thermal energy, $\mathrm{J}$ \\
	$R$ & Gas constant, $\mathrm{J/(mole\cdot K)}$ & $\real$ & Real numbers \\
	$\real_0^+$ & Non-negative real numbers & $R_{\text{chem}}$ & Arrhenius kinetic rate, $\mathrm{s^{-1}}$ \\
	$R_{\text{diff}}$ & Gas diffusion rate, $\mathrm{m^2/s}$ & $\mathbf{r}$ & Vector of deterministic variables used for transformation \\
	$S$ & Number of variate in PDD & $Sh$ & Sherwood number \\
	$T$ & Temperature, $\mathrm{K}$ & $\bar{\mathbf{T}}$ & Intermediate matrix for calculating $\breve{\bc}$, details in \cite{lee2024global} \\
	$T_{\text{avg}}$ & Average temperature at the pressure outlet, $\mathrm{K}$ & $\mathbf{W}$ & Diagonal matrix in cost function $\breve{\mathcal{K}}$ \\
	$w_1, w_2$ & Weighting factors in the objective function & $\mathbf{X}$ & Uncertain input vector \\
	$y(\cdot)$ & Generic output function & $\mathbf{Z}$ & Vector of transformed random variables \\
	$\mathbf{z}$ & Realization of $\bZ$ & & \\
	$\beta$ & Empirical unitless exponent & $\epsilon_{cp}$ & Maximum possible packing fraction \\
	$\epsilon_s$ & Empirical pressure constant & $\lambda$ & Weight in cost function $\breve{\mathcal{K}}$ \\
	$\mu_0, \sigma_0$ & Normalizing factors for the mean and standard deviation in the objective function & $\mathbf{\Phi}$ & Intermediate matrix for calculating $\breve{\bc}$, details in \cite{lee2024global} \\
	$\Psi(\cdot)$ & Univariate orthonormal polynomial function & $\tau$ & Collision stress, $\mathrm{Pa}$ \\
\end{longtable}
\end{appendices}
\bibliography{references}

\begin{thebibliography}{10}

\bibitem{acp-19-8523-2019}
M.~O. Andreae.
\newblock Emission of trace gases and aerosols from biomass burning -- an
  updated assessment.
\newblock {\em Atmospheric Chemistry and Physics}, 19(13):8523--8546, 2019.

\bibitem{beyer2007robust}
Hans-Georg Beyer and Bernhard Sendhoff.
\newblock Robust optimization--a comprehensive survey.
\newblock {\em Computer methods in applied mechanics and engineering},
  196(33-34):3190--3218, 2007.

\bibitem{boileau2008investigation}
Mathieu Boileau, St{\'e}phane Pascaud, Eleonore Riber, Benedicte Cuenot, LYM
  Gicquel, TJ~Poinsot, and M~Cazalens.
\newblock Investigation of two-fluid methods for large eddy simulation of spray
  combustion in gas turbines.
\newblock {\em Flow, Turbulence and Combustion}, 80:291--321, 2008.

\bibitem{burcat2005third}
Alexander Burcat and Branko Ruscic.
\newblock Third millenium ideal gas and condensed phase thermochemical database
  for combustion (with update from active thermochemical tables).
\newblock Technical report, Argonne National Lab.(ANL), Argonne, IL (United
  States), 2005.

\bibitem{chatterjee2019critical}
Tanmoy Chatterjee, Souvik Chakraborty, and Rajib Chowdhury.
\newblock A critical review of surrogate assisted robust design optimization.
\newblock {\em Archives of Computational Methods in Engineering}, 26:245--274,
  2019.

\bibitem{osti_1630426}
Mary~Ann Clarke and Jordan~M. Musser.
\newblock The {MFiX} particle-in-cell method ({MFiX-PIC}) theory guide.
\newblock Technical report, National Energy Technology Laboratory (NETL),
  Pittsburgh, PA, 2020.

\bibitem{cui2007fluidization}
Heping Cui and John~R Grace.
\newblock Fluidization of biomass particles: A review of experimental
  multiphase flow aspects.
\newblock {\em Chemical Engineering Science}, 62(1-2):45--55, 2007.

\bibitem{doucet2000sequential}
Arnaud Doucet, Simon Godsill, and Christophe Andrieu.
\newblock On sequential {Monte Carlo} sampling methods for {Bayesian}
  filtering.
\newblock {\em Statistics and Computing}, 10:197--208, 2000.

\bibitem{dryer1973high}
Frederick~L Dryer and I~Glassman.
\newblock High-temperature oxidation of {CO} and {CH4}.
\newblock In {\em Symposium (International) on Combustion}, volume~14, pages
  987--1003. Elsevier, 1973.

\bibitem{helton2003latin}
Jon~C Helton and Freddie~Joe Davis.
\newblock Latin hypercube sampling and the propagation of uncertainty in
  analyses of complex systems.
\newblock {\em Reliability Engineering \& System Safety}, 81(1):23--69, 2003.

\bibitem{huang2006robust}
Beiqing Huang and Xiaoping Du.
\newblock A robust design method using variable transformation and
  {Gauss--Hermite} integration.
\newblock {\em International Journal for Numerical Methods in Engineering},
  66(12):1841--1858, 2006.

\bibitem{huang2007analytical}
Beiqing Huang and Xiaoping Du.
\newblock Analytical robustness assessment for robust design.
\newblock {\em Structural and multidisciplinary optimization}, 34:123--137,
  2007.

\bibitem{janssen2013monte}
Hans Janssen.
\newblock {Monte-Carlo} based uncertainty analysis: Sampling efficiency and
  sampling convergence.
\newblock {\em Reliability Engineering \& System Safety}, 109:123--132, 2013.

\bibitem{jin2003use}
Ruichen Jin, Xiaoping Du, and Wei Chen.
\newblock The use of metamodeling techniques for optimization under
  uncertainty.
\newblock {\em Structural and Multidisciplinary Optimization}, 25:99--116,
  2003.

\bibitem{KORKERD2024119262}
Krittin Korkerd, Zongyan Zhou, Ruiping Zou, Pornpote Piumsomboon, and Benjapon
  Chalermsinsuwan.
\newblock Numerical investigation of mixing and heat transfer of mixed biomass
  and silica sand particles in a bubbling fluidized bed combustor.
\newblock {\em Powder Technology}, 433:119262, 2024.

\bibitem{LEE2022103218}
Dongjin Lee, Ramin Jahanbin, and Sharif Rahman.
\newblock Robust design optimization by spline dimensional decomposition.
\newblock {\em Probabilistic Engineering Mechanics}, 68:103218, 2022.

\bibitem{lee2023bi}
Dongjin Lee and Boris Kramer.
\newblock Bi-fidelity conditional value-at-risk estimation by dimensionally
  decomposed generalized polynomial chaos expansion.
\newblock {\em Structural and Multidisciplinary Optimization}, 66(2):33, 2023.

\bibitem{lee2023multifidelity}
Dongjin Lee and Boris Kramer.
\newblock Multifidelity conditional value-at-risk estimation by dimensionally
  decomposed generalized polynomial chaos-{Kriging}.
\newblock {\em Reliability Engineering \& System Safety}, 235:109208, 2023.

\bibitem{lee2024global}
Dongjin Lee, Elle Lavichant, and Boris Kramer.
\newblock Global sensitivity analysis with limited data via sparsity-promoting
  {D-MORPH} regression: Application to char combustion.
\newblock {\em Journal of Computational Physics}, 511:113116, 2024.

\bibitem{lee2021robust}
Dongjin Lee and Sharif Rahman.
\newblock Robust design optimization under dependent random variables by a
  generalized polynomial chaos expansion.
\newblock {\em Structural and Multidisciplinary Optimization},
  63(5):2425--2457, 2021.

\bibitem{lee2023high}
Dongjin Lee and Sharif Rahman.
\newblock High-dimensional stochastic design optimization under dependent
  random variables by a dimensionally decomposed generalized polynomial chaos
  expansion.
\newblock {\em International Journal for Uncertainty Quantification}, 13(4),
  2023.

\bibitem{mckinnon1998convergence}
Ken~IM McKinnon.
\newblock Convergence of the {Nelder}--{Mead} simplex method to a nonstationary
  point.
\newblock {\em SIAM Journal on Optimization}, 9(1):148--158, 1998.

\bibitem{noh2009reliability}
Yoojeong Noh, KK~Choi, and Liu Du.
\newblock Reliability-based design optimization of problems with correlated
  input variables using a {Gaussian Copula}.
\newblock {\em Structural and Multidisciplinary Optimization}, 38:1--16, 2009.

\bibitem{osman2023cost}
Ahmed~I Osman, Lin Chen, Mingyu Yang, Goodluck Msigwa, Mohamed Farghali, Samer
  Fawzy, David~W Rooney, and Pow-Seng Yap.
\newblock Cost, environmental impact, and resilience of renewable energy under
  a changing climate: a review.
\newblock {\em Environmental Chemistry Letters}, 21(2):741--764, 2023.

\bibitem{rahman2008polynomial}
Sharif Rahman.
\newblock A polynomial dimensional decomposition for stochastic computing.
\newblock {\em International Journal for Numerical Methods in Engineering},
  76(13):2091--2116, 2008.

\bibitem{rahman2009extended}
Sharif Rahman.
\newblock Extended polynomial dimensional decomposition for arbitrary
  probability distributions.
\newblock {\em Journal of Engineering Mechanics}, 135(12):1439--1451, 2009.

\bibitem{ren2013robust}
Xuchun Ren and Sharif Rahman.
\newblock Robust design optimization by polynomial dimensional decomposition.
\newblock {\em Structural and Multidisciplinary Optimization}, 48:127--148,
  2013.

\bibitem{ren2016reliability}
Xuchun Ren, Vaibhav Yadav, and Sharif Rahman.
\newblock Reliability-based design optimization by adaptive-sparse polynomial
  dimensional decomposition.
\newblock {\em Structural and Multidisciplinary Optimization}, 53:425--452,
  2016.

\bibitem{robert2013monte}
Christian Robert and George Casella.
\newblock {\em {Monte} Carlo Statistical Methods}.
\newblock Springer Science \& Business Media, 2013.

\bibitem{shen2016polynomial}
Dongying~E Shen and Richard~D Braatz.
\newblock Polynomial chaos-based robust design of systems with probabilistic
  uncertainties.
\newblock {\em AIChE Journal}, 62(9):3310--3318, 2016.

\bibitem{snider2001incompressible}
Dale~M Snider.
\newblock An incompressible three-dimensional multiphase particle-in-cell model
  for dense particle flows.
\newblock {\em Journal of Computational Physics}, 170(2):523--549, 2001.

\bibitem{steiner2024particulate}
Dominik Steiner and Christof Lanzerstorfer.
\newblock Particulate emissions from biomass power plants: a practical review
  and measurement uncertainty issues.
\newblock {\em Clean Technologies and Environmental Policy}, 26(4):1039--1048,
  2024.

\bibitem{sudret2008global}
Bruno Sudret.
\newblock Global sensitivity analysis using polynomial chaos expansions.
\newblock {\em Reliability Engineering \& System Safety}, 93(7):964--979, 2008.

\bibitem{tu1999new}
Jian Tu, Kyung~K Choi, and Young~H Park.
\newblock A new study on reliability-based design optimization.
\newblock {\em Journal of Mechanical Design}, 121(4):557--564, 12 1999.

\bibitem{wang2021efficient}
Jian Wang, Zhili Sun, and Runan Cao.
\newblock An efficient and robust {K}riging-based method for system reliability
  analysis.
\newblock {\em Reliability Engineering \& System Safety}, 216:107953, 2021.

\bibitem{weickum2006multi}
Gary Weickum, Mike Eldred, and Kurt Maute.
\newblock Multi-point extended reduced order modeling for design optimization
  and uncertainty analysis.
\newblock In {\em 47th AIAA/ASME/ASCE/AHS/ASC Structures, Structural Dynamics,
  and Materials Conference 14th AIAA/ASME/AHS Adaptive Structures Conference
  7th}, page 2145, 2006.

\bibitem{xie2021study}
Jun Xie, Wenqi Zhong, and Yingjuan Shao.
\newblock Study on the char combustion in a fluidized bed by {CFD-DEM}
  simulations: Influences of fuel properties.
\newblock {\em Powder Technology}, 394:20--34, 2021.

\bibitem{xie2019coupling}
Jun Xie, Wenqi Zhong, Yingjuan Shao, and Kaixi Li.
\newblock Coupling of {CFD-DEM} and reaction model for 3d fluidized beds.
\newblock {\em Powder Technology}, 353:72--83, 2019.

\bibitem{yang2021hybrid}
Meide Yang, Dequan Zhang, Chao Jiang, Xu~Han, and Qing Li.
\newblock A hybrid adaptive {K}riging-based single loop approach for complex
  reliability-based design optimization problems.
\newblock {\em Reliability Engineering \& System Safety}, 215:107736, 2021.

\bibitem{youn2005performance}
Byeng~D Youn, Kyung~K Choi, and Kiyoung Yi.
\newblock Performance moment integration {(PMI)} method for quality assessment
  in reliability-based robust design optimization.
\newblock {\em Mechanics Based Design of Structures and Machines},
  33(2):185--213, 2005.

\bibitem{zou2006direct}
T~Zou and Sankaran Mahadevan.
\newblock A direct decoupling approach for efficient reliability-based design
  optimization.
\newblock {\em Structural and Multidisciplinary Optimization}, 31:190--200,
  2006.

\end{thebibliography}
\bibliographystyle{plain}

%
\end{document}